\documentclass{article}
\usepackage{amsthm, amsmath,amssymb, amscd}
\usepackage{txfonts}
\usepackage[dvips]{graphicx, color}
\usepackage{latexsym,amsxtra,amsfonts,verbatim}
\usepackage{lipsum}     % dummy text
\usepackage[makeroom]{cancel}

\usepackage[dvips]{graphicx, color}
\usepackage[title]{appendix}

\newcommand{\xdownarrow}[1]{%
  {\left\downarrow\vbox to #1{}\right.\kern-\nulldelimiterspace}
}
\usepackage[new]{old-arrows}

\usepackage{bbm}

\theoremstyle{plain}
\newtheorem{thm}{Theorem}[section]
\newtheorem{lem}[thm]{Lemma}
\newtheorem{sublem}[thm]{Sublemma}
\newtheorem{prop}[thm]{Proposition}
\newtheorem{cor}[thm]{Corollary}

\theoremstyle{definition}

\newtheorem{conj}{Conjecture}[section]
\newtheorem{exam}{Example}[section]

\theoremstyle{remark}

\newcommand{\cE}{{\mathcal E}}

\newcommand{\C}{{\Bbb C}}

\newcommand{\F}{{\Bbb F}}

\renewcommand{\O}{{\Bbb O}}

\newcommand{\Q}{{\Bbb Q}\hspace{.06em}}
\newcommand{\R}{{\Bbb R}}

\newcommand{\Z}{{\Bbb Z}}
\def\O{\mathcal O}

\def\={\:=\:}  \def\+{\,+\,}

\def\a{\alpha} \def\b{{\beta}}  \def\ba{\overline\a}

\def\be{\begin{equation}}   \def\ee{\end{equation}}
\def\bes{\begin{equation*}}   \def\ees{\end{equation*}}
\def\ba{\begin{aligned}}   \def\ea{\end{aligned}}
\def\bc{\begin{cases}}   \def\ec{\end{cases}}
\def\bp{\begin{proof}}   \def\ep{\end{proof}}

\newcommand{\Res}{\mathrm{Res}}

\newcommand{\Aut}{\mathrm{Aut}}

\def\SL{\mathrm{SL}}

\newcommand{\ov}{\overline}
\def\qqan{\qquad\mathrm{and}\qquad}
\def\qan{\quad\mathrm{and}\quad}

\newcommand{\dis}{\displaystyle}

\def\smm{\smallsetminus}
\def\s{\sigma}

\def\ga{\gamma}

\def\SL{\mathrm{SL}}

\def\ov{\overline}
\def\De{\Delta}

\def\bbm1{\mathbbm 1}

\def\De{\Delta}

\def\wh{\widehat}
\def\be{\begin{equation}}   \def\ee{\end{equation}}
\def\bes{\begin{equation*}}   \def\ees{\end{equation*}}
\def\bea{\begin{equation}\begin{aligned}}   
\def\eea{\end{aligned}\end{equation}}

\def\sm{\smm}

\def\tot{\mathrm{tot}}

\def\Pic{\mathrm{Pic}}

\def\om{\omega}

\def\tot{\mathrm{tot}}

\def\bm{\begin{matrix}}
\def\em{\end{matrix}}
\def\bpm{\begin{pmatrix}}
\def\epm{\end{pmatrix}}
\def\Hom{\mathrm{Hom}}

\def\bl{\big(}
\def\br{\big)}

\def\de{\delta}

\def\Z{\mathbb Z} \def\R{\mathbb R}    
\def\={\;=\;}  \def\+{\,+\,}  \def\C{\Bbb C}    \def\Q{\Bbb Q}  
\def\Z{\Bbb Z}  \def\F{\Bbb F}

\def\aA{\mathcal A}  \def\bB{\mathcal B}

\begin{document}

\title{\bf  Riemann Hypothesis\\ for Non-Abelian Zeta Functions\\
 of Curves over Finite Fields} 
\author{Lin WENG}  
\date{}
\maketitle
\begin{abstract} In this paper, we develop some basic techniques towards the Riemann hypothesis for higher rank non-abelian zeta functions of  an integral  regular projective curve of genus $g$ over a finite field $\F_q$. As an application of the Riemann hypothesis for these genuine zeta functions, we obtain some explicit bounds on the fundamental non-abelian $\a$- and $\b$-invariants of $X/\F_q$ in terms of $X$ and $n,\, q$ and $g$:
$$\a_{X,\F_q;n}(mn)=\sum_{V}\frac{q^{h^0(X,V)}-1}{\#\Aut(V)}\qqan \b_{X,\F_q;n}(mn):=\sum_{V}\frac{1}{\#\Aut(V)}\qquad(m\in \Z)$$
where $V$ runs through all rank $n$ semi-stable $\F_q$-rational vector bundles on $X$ of degree $mn$. In particular,
$$
\prod_{k=1}^{n}\frac{\ \bl \sqrt q^k-1\br^{2g-1}\ }{(\sqrt q^k+1)}\leq q^{-\binom{n}{2}(g-1)} \b_{X,\F_q;n}(0)
\leq \prod_{k=1}^{n}\frac{\ \bl 1+\sqrt q^k\br^{2g-1}\ }{(\sqrt q^k-1)},
$$
Finally, we demonstrate that the related bounds in lower ranks in turn  play a central role in establishing the Riemann hypothesis for higher rank zetas.
\end{abstract}

\tableofcontents

\section{Special uniformity of zetas}
Special uniformity for zeta functions of curves over finite fields is conjectured in \cite{W} and established in \cite{WZ2}, with the help of the result in \cite{MR}. In this section, we recall some basic constructions involved.

\subsection{Non-abelian zeta function of a curve over a finite field}

For a fixed positive integer $n\geq 1$, the {\it rank $n$ non-abelian zeta function} of a projective regular integral curve $X$ over $\F_q$ is defined in \cite{W}\footnote{While this paper is fundamental to the field, it has never been submitted for a formal publication.} as
\be
\wh\zeta_{X,\F_q;n}(s)=\sum_{m=0}^\infty\sum_{V}\frac{q^{h^0(X,\cE)-1}}{\#\Aut\, \cE}\bl q^{-s}\br^{\chi(X,\cE)}\hskip 2.0cm\bl\Re(s)>1\br
\ee
where $\cE$ (in the second summation) runs through rank $n$ semi-stable vector bundles on $X/\F_q$ of degree $mn$. This definition is a modification of the previous one in \cite{W0} in which $\cE$ is allowed to run over all  rank $n$ semi-stable vector bundles on $X/\F_q$ of degree $d\ (d\geq 0)$. However, despite the fact that this original definition in \cite{W0} would yield a rational function satisfying the standard functional equation, it fails to satisfy any form of the Riemann hypothesis. To remedy this, motivated by Drinfeld's work \cite{D} on counting two-dimensional irreducible representations of the fundamental group of a curve over a finite field, we introduce a restriction on the degrees of $\cE$, that is, the degrees of $V$ are required to be divided by $n$, the rank of $\cE$. By using the Riemann-Roch theorem, cohomology duality and  vanishing theorem for semi-stable bundles on curves, tautologically, we have
\begin{thm}[$\zeta$ Properties \cite{W}]  The rank $n$-zeta function $\wh\zeta_{X,\F_q;n}(s)$ of a genus $g$  integral, regular projective curve $X$ over $\F_q$ satisfies the following canonical zeta properties:
 \begin{enumerate}
\item [(0)]  $\wh\zeta_{X,\F_q;1}(s)$ coincides with the (completed) Artin zeta $\wh\zeta_{X/\F_q}(s)$ of $X/\F_q$. 
\item [(1)] $\wh\zeta_{X,\F_q;n}(s)$ is a rational function in $T:=\bl q^{-s}\br^n$.
\item[(2)] {\rm (Functional\ equation)}  $\wh\zeta_{X,\F_q;1}(s)$ satisfies the standard functional equation
\be
\wh\zeta_{X,\F_q;n}(1-s)=\wh\zeta_{X,\F_q;n}(s).
\ee
\item[(3)] {\rm (Singularities)} As a rational function of $T$, $\wh\zeta_{X,\F_q;n}(s)$ admits only two singularities, namely, two simple poles at $T=0,1$. Furthermore, the residue at $T=1$ admits the following geo-arithmetic interpretation:
\be
\Res^{~}_{T=1}\wh\zeta_{X,\F_q;n}(s)=\b_{X,\F_q;n}(0):=\sum_{\cE}\frac{1}{\#\,\Aut(\cE)}
\ee
where $\cE$ runs over rank $n$ semi-stable vector bundles on $X/\F_q$ of degree 0.
\end{enumerate}
\end{thm}

Indeed, (0) can be deduced by expressing the Artin zeta $\wh\zeta_{X/\F_q}(s)$ of $X/\F_q$ as a sum on the rationally equivalence classes of divisors, or better, the rational line bundles, of 
non-negative degrees on $X/\F_q$. Furthermore,
 if we introduce the non-abelian geo-arithmetic $\a$- and $\b$-invariants of the curve $X$ over $\F_q$ associated to rank $n$ semi-stable vector bundles  by
 \be\a_{X,\F_q;n}(d)=\sum_{\cE}\frac{q^{h^0(X,\cE)}-1}{\#\Aut(\cE)}\qqan \b_{X,\F_q;n}(d):=\sum_{\cE}\frac{1}{\#\Aut(\cE)}\qquad(\forall d\in \Z)\ee
 where $\cE$ rums through rank $n$ semi-stable $\F_q$-rational vector bundles of degree $d$ on $X$,\footnote{The beta invariant was first introduced in \cite{HN}.}
then by the vanishing theorem for semi-stable vector bundles $\cE$ and the Riemann-Roch theorem,  we conclude that
 \be
 \a_{X,\F_q;n}(d)=0\quad(d<0)\qqan \a_{X,\F_q;n}(d)=\b_{X,\F_q;n}(d)\cdot \Big(q^{d-n(g-1)}-1\Big)\quad (\forall d\geq 2n(g-1)).
 \ee
 In addition, directly from the definition, we have
 \be
  \b_{X,\F_q;n}(mn)=\b_{X,\F_q;n}(0)\qquad(\forall m\in\Z).
 \ee
Therefore,  by the cohomology duality and  standard $\zeta$-technique for curves, we have the following
 \begin{thm}[\cite{W}]\label{thm1.2} The rank $n$-zeta function $\wh\zeta_{X,\F_q;n}(s)$ of a genus $g$ projective regular integral curve $X$ over $\F_q$ is given by, with $Q=q^n$,
 $$\ba
 &\wh Z_{X,\F_q;n}(T):=\wh \zeta_{X,\F_q;n}(s)\\
 =&
 \sum_{m=0}^{g-2}\a_{X,n}(mn)\left(T^{m-(g-1)}+Q^{(g-1)-m}T^{(g-1)-m}\right) +\a_{X,n}\bl n(g-1)\br+\frac{(Q-1)\b_{X,n}(0)\cdot T}{(1-T)(1-QT)}.
 \ea
 $$
 \end{thm}
Here, of course, when $g-2<0$, resp. $g-1<0$, the first sum, resp. the first sum and the second term, is understood as 0.
 In particular, when $n=1$, we have recovered the following basic, but less well-known geo-arithmetic interpretation of the Artin Zeta function of the curve $X/\F_q$, 
 $$\ba
 &\wh Z_{X/\F_q}(t):=\wh \zeta_{X/\F_q}(s)\\
 =&
 \sum_{m=0}^{g-2}\a_{X/\F_q}(m)\left(t^{m-(g-1)}+q^{(g-1)-m}t^{(g-1)-m}\right) +\a_{X/\F_q}\bl g-1\br+\frac{(q-1)\b_{X/\F_q}(0)\, t}{(1-t)(1-qt)}.
 \ea
 $$
 where, $t=q^{-s}$, and to simplify our notation, we have set $\a_{X/\F_q}(d):=\a_{X,\F_q;1}(d)$ and $\b_{X/\F_q}(d):=\b_{X,\F_q;1}(d)$.
 
This theorem clearly implies the zeta properties on rationality, the functional equation and the singularities in the previous theorem. Furthermore, if we set
 \be\label{eq5}
 \wh Z_{X,\F_q;n}(T)=\frac{P_{X,\F_q;n}(T)}{(1-T)(1-QT)\cdot T^{g-1}}
 \ee
Then $P_{X,\F_q;n}(T)$ is a degree $2g$ polynomial in $T$ with real coefficients whose  leading coefficient and constant term are given by $Q^{g}\a_{X,\F_q;n}(0)$ and  $\a_{X,\F_q;n}(0)$, respectively.

After examining many examples in lower ranks, in \cite{W}, we formulate the following
\begin{conj}[Riemann Hypothesis] {\it The rank $n$-zeta function $\wh\zeta_{X,\F_q;n}(s)$ of a projective regular integral curve $X$ over $\F_q$ satisfies the Riemann hypothesis. That is,
all  roots of $p^{~}_{X,\F_q;n}(s):=P_{X,\F_q;n}(q^{-s})$ lies on the line $\Re(s)=\frac{1}{2}$.}
\end{conj}
Obviously, this is equivalent to the condition that all reciprocal roots of $P_{X,\F_q;n}(T)$ are of norm $Q^{\frac{1}{2}}$. Still, there is an apparently weak but equivalent  form which claims that all reciprocal roots of $P_{X,\F_q;n}(T)$ are not real, thanks to the functional equation of the non-abelian zeta $\wh\zeta_{X,\F_q;n}(s)$. 

The first break-through  in this direction is the following result on elliptic curves of Zagier and myself \cite{WZ1}, relying on some basic properties of Atiyah bundles \cite{A} and a heavy use of combinatorics:
\begin{thm}[\cite{WZ1}] Let $E$ be an elliptic curve  over $\F_q$. Then the rank $n$ zeta function
$\wh\zeta_{E,\F_q;n}(s)$ of $E$ satisfies the Riemann hypothesis.
\end{thm}

\subsection{$\SL_n$-zeta functions of a curve over a finite field}

Let $X$ be an integral regular projective curve over $\F_q$ and let $G$ be a  split connected reductive algebraic group of rank $r$ over  $\F_q(X)$,  the function field of $X/\F_q$. Let 
\be
\Big(V,\langle\cdot,\cdot\rangle,\Phi=\Phi^+\cup\Phi^-,\Delta=\{\alpha_1,\dots,\alpha_r\},
\varpi:=\{\varpi_1,\dots,\varpi_r\},W\Big)
\ee
be the  {\it root system} associated to  a fixed split minimal parabolic subgroup~$P_0$ of $G$ and its maximal split torus~$T$. Here, as usual, $V$ can be identified with the real vector space of
$\mathbb R$-span of rational characters of $T$, and is equipped with a natural inner product $\langle\cdot,\cdot\rangle$, with which we may and hence will identify $V$ with its dual $V^*:=\Hom_\R(V,\R)$. In addition,  $\Phi^+\subset V$, resp. $\Phi^-:=-\Phi^+$ , denotes the set of the so-called {\it positive roots},  resp. {\it negative roots}, 
$\Delta\subset V$, resp. $\varpi\subset V$, denotes the set of {\it simple roots}, resp. of {\it fundamental weights}, 
and $W$ denotes the  {\it  Weyl group} generated by the reflections $\s_\a$ $(\a\in\De)$.  By definition, the fundamental weights are characterized by 
the formula 
\be
\langle\varpi_i,\alpha_j^\vee\rangle=\delta_{ij}\qquad(\forall i,j=1,2,\dots,r),
\ee 
where, for each root $\alpha\in\Phi$,  $\alpha^\vee:=\frac{2}{\langle\alpha,\alpha\rangle}\,\alpha$ denotes the corresponding  {\it coroot}.  We also define the {\it Weyl vector} $\rho$ by  
\be
\rho=\frac{1}{2}\sum_{\alpha\in\Phi^+}\alpha,
\ee 
and introduce a {\it coordinate system} on~$V$ (with respect to the base $\{\varpi_1,\dots,\varpi_r\}$ of $V$ and the vector $\rho$) by writing 
an element  $\lambda \in V$ in the form
\be
\lambda\=\sum_{j=1}^r(1-s_j)\varpi_j\=\rho-\sum_{j=1}^rs_j\varpi_j\,.
\ee
This, in turn, induces natural  identifications of $V$ and $V_\C=V\otimes_\R\C$ with $\R^r$ and $\C^r$, respectively. 
For each Weyl element $w\in W$, we set 
\be
\Phi_w:=\Phi^+\,\cap\, w^{-1}\Phi^-,
\ee 
be the collection of positive roots whose $w$-images are negative. It is well-known that the cardinality of $\Phi_w$ coincides with the {\it length} $\ell_w$ of $w$, i.e. the minimal number expressing $w$ in terms of products of $\s_\a\ (\a\in\De)$.

\medskip
As usual, by a {\it standard parabolic subgroup} of $G$, we mean a parabolic subgroup of $G$ that contains the fixed minimal parabolic subgroup $P_0$. From Lie theory (see e.g.,~\cite{Hu}), there is an one-to-one 
correspondence between standard parabolic subgroups $P$ of $G$ and subsets $\Delta_P$ 
of $\Delta$. In particular, if $P$ is maximal, we may and will write 
$\Delta_P=\Delta\sm \{\alpha_p\}$ for a certain unique $p=p(P)\in\{1,\dots,r\}$. For such a 
standard parabolic subgroup $P$, denote by $V_P$  the $\mathbb R$-span of rational 
characters of the maximal split torus $T_P$ contained in  $P$, by $V_P^*$ its dual space, 
and by $\Phi_P\subset V_P$ the set of non-trivial characters of $T_P$ occurring in the space 
$V$. Then, by the  standard theory of reductive groups (see e.g.,~\cite{B}), $V_P$, resp. $V_P^*$, admits a 
canonical embedding in $V$, resp. in $V^*$, which is 
known to be orthogonal to the fundamental weight $\varpi_p$, and hence $\Phi_P$ can be viewed 
as a subset of $\Phi$. Set, accordingly, 
\be
\Phi^+_P:=\Phi^+\,\cap\,\Phi_P,\qquad\rho_P:=\frac{1}{2}\sum_{\alpha\in\Phi_P^+}\alpha\qqan c_P:=2\langle\varpi_p-\rho_P,\alpha_p^\vee\rangle.
\ee 

\medskip
Now, for an integral  regular projective curve $X$ of genus $g$ over a finite field $\mathbb F_q$, 
in \cite{W3}, motivated by the study of zeta functions for number fields,\footnote{For number 
fields, the analogue of the two functions to be  introduced below are special kinds of Eisenstein 
periods, defined as integrals of Eisenstein series over moduli spaces of semi-stable lattices. 
For details, see~\cite{W2}.}  for a connected split reductive algebraic group~$G$, and its maximal standard parabolic subgroup $P$ (defined over the function field of $X/\F_q$), we introduce the 
{\it period of $G$} and  the {\it period of $(G,P)$} for $X/\F_q$ by  
$$
\omega_{X,\F_q}^G(\lambda):=\sum_{w\in W}\frac{1}{\prod_{\alpha\in\Delta}
(1-q^{-\langle w\lambda-\rho,\alpha^\vee\rangle})}
\prod_{\alpha\in\Phi_w}\frac{\wh\zeta_{X/\F_q}(\langle\lambda,\alpha^\vee\rangle)}
{\wh\zeta_{X/\F_q}(\langle\lambda,\alpha^\vee\rangle+1)} 
$$
and  
$$
\begin{aligned}
\omega^{G,P}_{X/\F_q}(s)
\;&:=\; \mathrm{Res}_{\langle\lambda-\rho,\,\alpha^\vee\rangle=0,\;\alpha\in\Delta_P}
\omega_{X,\F_q}^G(\lambda)\bigr|_{s_p=s} \\
&\phantom{:}=\;\mathrm{Res}_{s_r=0}\cdots \mathrm{Res}_{s_{p+1}=0}\mathrm{Res}_{s_{p-1}=0}\cdots
\mathrm{Res}_{s_1=0}\,\omega_{X/\F_q}^G(\lambda)\bigr|_{s_p=s} \;,
\end{aligned}
$$ 
respectively,
where $s$ is a complex variable\footnote{We should warn the reader that in~\cite{W3}, \cite{W2}
and~\cite{W} a different normalization is used, with the argument of $\omega_{X,\F_q}^{G,P}$ (and
later of $\zeta_X^{G,P}$) being given by $s=c_p(s_p-1)$ ($\,=n(s_p-1)$ in the special case
$(G,P)=(\SL_n,P_{n-1,1})$) rather than $s=s_p$ as chosen here.  With the normalization used here 
 the functional equation relates $s$~and $1-s$ rather than $s$~and~$-n-s$.}
and where for the last equality we have used the facts that
\be
\langle\rho,\alpha^\vee\rangle=1\qquad(\forall \alpha\in\Delta)\qqan\langle\varpi_i,\alpha_j^\vee\rangle=\delta_{ij}\quad(\forall 1\leq i,\,j\leq r).\ee
As proved in \cite{Ko, W3}, the ordering of taking residues along singular hyperplanes 
$\langle\lambda-\rho,\alpha^\vee\rangle=0$ for $\alpha\in\Delta_P$ does not affect the 
outcome, so that the definition is independent of the numbering of the simple roots used in the definition.

\medskip
To get the zeta function associated to $(G,P)$ for $X/\F_q$, certain normalizations should be made.
For this purpose, write $\displaystyle{\omega_{X,\F_q}^G(\lambda)\={\sum}_{w\in W}T_w(\lambda)}$, 
where, for each $w\in W$,
$$
T_w(\lambda):=\frac{1}{\prod_{\alpha\in\Delta}(1-q^{-\langle w\lambda-\rho,\alpha^\vee\rangle})}
\prod_{\alpha\in\Phi_w}\frac{\wh\zeta_{X/\F_q}(\langle\lambda,\alpha^\vee\rangle)}
{\wh\zeta_{X/\F_q}(\langle\lambda,\alpha^\vee\rangle+1)}.
$$
Accordingly, we need to undertand the residue
$$\mathrm{Res}_{\langle\lambda-\rho,\,\alpha^\vee\rangle=0,\;\alpha\in\Delta_P}T_w(\lambda).$$
Clearly, we care only about those elements~$w\in W$ (which we will call {\it special}) that give
non-trivial residues, namely, those satisfying the condition that
$\mathrm{Res}_{\langle\lambda-\rho,\,\alpha^\vee\rangle=0,\;\alpha\in\Delta_P}
T_w(\lambda)\not\equiv 0.$ This can happen only if all singular hyperplanes  are of one of the 
following two forms:
\begin{enumerate}
\renewcommand{\labelenumi}{(\arabic{enumi})}
\item $\langle w\lambda-\rho,\alpha^\vee\rangle=0$ for some $\alpha\in\Delta$, giving a 
simple pole of the rational factor 
$\frac{1}{\prod_{\alpha\in\Delta}(1-q^{-\langle w\lambda-\rho,\alpha^\vee\rangle})}$;
\item  $\langle \lambda,\alpha^\vee\rangle=1$ for some $\alpha\in\Phi_w$, giving a simple 
pole of the zeta factor $\wh\zeta_{X/\F_q}(\langle\lambda,\alpha^\vee\rangle)$ appeared in the numerator within the second group of factors, the product.
\end{enumerate}

\noindent
For special $w\in W$, and $(k,h)\in\mathbb Z^2$, following~\cite{Ko}
(see also~\cite{W3}) we define 
\begin{align} N_{P,w}(k,h)\,:=\,&\#\Big\{\alpha\in w^{-1}\Phi^-\,:\,
   \langle\varpi_p,\alpha^\vee\rangle=k,\;\langle\rho,\alpha^\vee\rangle=h\Big\} \notag\\[0.3em]
M_P(k,h):=&\max_{\text{$w$ special}}\Bigl\{N_{P,w}(k,h-1)-N_{P,w}(k,h)\Bigr\}\,.\notag\\
=\,&N_{P,w_0}(k,h-1)-N_{P,w_0}(k,h) \,,\label{Cor87} \end{align}
where $w_0$ is the longest element of the Weyl group $W$. Indeed,  the last equality is guaranteed by Corollary~8.7 of~\cite{KKS}.
Note that $M_P(k,h)=0$ for almost all but finitely many pairs of integers $(k,h)$, so it makes 
sense to introduce the product 
\be\label{defDGP}
D^{G,P}_{X,\F_q}(s):=\prod_{k=0}^\infty\prod_{h=2}^\infty\wh\zeta_{X,\F_q}(kn(s-1)+h)^{M_P(k,h)}.
\ee
Following~\cite{W2,W3}, we define the {\it zeta function of $X/\F_q$ associated to $(G,P)$} by 
\begin{equation}
{\wh\zeta}_{X,\F_q}^{G,P}(s):=q^{(g-1)\dim N_u(P_0)}\cdot D_{X,\F_q}^{G,P}(s)\cdot \omega_{X,\F_q}^{G,P}(s)\,. 
\label{eq: UE}
\end{equation}
Here $N_u(P_0)$ denote the nilpotent radical of the minimal parabolic subgroup $P_0$ of $G$.
\vskip 0.20cm
\noindent
{\it Remark.} For special $w\in W$, even after taking residues, there are some zeta factors  $\wh\zeta_{X/\F_q}(ks+h)$ left
in the denominator of $\mathrm{Res}_{\langle\lambda-\rho,\,\alpha^\vee\rangle=0,\;\alpha\in\Delta_P}T_w(\lambda)$. 
The reason for introducing the factor $D^{G,P}_X(s)$ in our normalization of the zeta functions, based on formulas in~\cite{Ko} 
and~\cite{W3}, is to clear up all of those zeta factors (appearing in the denominators) associated to special Weyl elements.

In particular, we have the following
\begin{thm}[Functional Equation\cite{W}] For an integral  regular projective curve $X$ over $\F_q$,
\be
{\wh\zeta}_{X,\F_q}^{G,P}(c_P-s)={\wh\zeta}_{X,\F_q}^{G,P}(s).
\ee
\end{thm}
The proof follows \cite{Ko} closely, where the functional equation is established for the parallel  structures on the so-called $(G,P)$-zeta function of number fields $F$.

With all these, we are now ready to introduce  the $\SL_n$-zeta function of $X/\F_q$
by specializing to the case when $G$ is the special linear group $\SL_n$ 
and $P$ is the maximal parabolic subgroup $P_{n-1,1}$ consisting of matrices whose final row 
vanishes except for its last entry, corresponding to the ordered partition $(n-1)+1$ of~$n$. 
That is to say, the {\it $\SL_n$-zeta function ${\wh\zeta}_{X,\F_q}^{\SL_n}(s)$ of $X/\F_q$} is defined to be 
\begin{equation}
{\wh\zeta}_{X,\F_q}^{\SL_n}(s):=\;{\wh\zeta}_{X,\F_q}^{\SL_n,\,P_{n-1,1}}(s):= q^{\frac{n(n-1)}{2}(g-1)} \,\cdot\, D^{\SL_n,P_{n-1,1}}(s) \,\cdot\, \omega_{X}^{(\SL_n,P_{n-1,1})}(s)\;.
\label{eq: UG}
\end{equation}
 As the first step to understand this zeta function, we have the following

\begin{lem}[Lemma 5 of \cite{WZ2}] The finction $D^{\SL_n,P_{n-1,1}}(s)$ is given by
\be \label{DSLn}
D^{\SL_n,P_{n-1,1}}(s) \= \prod_{k=2}^{n-1}\wh\zeta_{X/\F_q}(k) \,\cdot\, \wh\zeta_{X/\F_q}(ns).
\ee
\end{lem}
Motivated by our study on the parallel structures for number fields, after verifying some concrete examples, in \cite{W}, we formulate the following
\begin{conj}[Special Uniformity of Zetas] For an integral regular projective curve $X$ on $\F_q$ of genus $g$,
up to some constant factor depending only on $n$ and $g$, we have
\be
\wh\zeta_{X,\F_q;n}(s)=\wh\zeta_{X/\F_q}^{\SL_n}(s)
\ee
\end{conj}

For number fields $F$, this uniformity of zeta functions is established using Mellin transforms to write down the rank $n$ non-abelian zeta function of $F$ in terms of integrations of Epstein zeta functions over the moduli space of semi-stable $\O_F$-lattices of rank $n$ and degree zero. But Epstein zeta function is a special kind of Eisenstein series, which can be realized as the residue
of the Siegel-Langlands Eisenstein series associated to the constant function on the Levi subgroup of the minimal parabolic subgroup $P_{1,1,\ldots,1}$ corresponding to the decomposition $n=\overbrace{1+\ldots +1}^{n\ times}$. Furthermore, the moduli space of of semi-stable $\O_F$-lattices of rank $n$ and degree zero can be identified with the truncated domain of Arthur type within the fundamental domain of $\SL_n(\Z)$, a deep structural result parallel to that of Lafforgue for function fields \cite{La}. Consequently, with an use of relative trace formula yields the desired zeta uniformity. For details, please refer to Chapter 15 of \cite{W2}.

To pave the same path to establish the special uniformity of zetas for function fields, the first difficulty is that the analogue construction of Mellin transform has yet to be developed (see however a work of K. Adachi at Kyushu university on \lq\lq Rankin-Selberg \& Zagier Methods for Function Fields over Finite Fields").

\subsection{Special uniformity of zeta functions}

As said, the special uniformity of zetas claims that, for a global field $F$, the geometrically defined rank $n$ zeta function $\wh\zeta_{F,n}(s)$  coincides with the Lie theoretically defined 
$\SL_n$-zeta function $\wh\zeta_{F}^{\SL_n}(s)$.
 When $F$ is a number field, this conjectured in confirmed in \cite{W2} using the theories of Eisenstein series of Siegel (classical) and Langlands (modern), Arthur's analytic truncation and geo-arithmetic truncation of stability, and relative trace formula. When $F$ is a function field, a totally different approach has been used, thanks to an unexpected work of Mozgovoy-Reineke \cite{MR}. The uniformity of zetas for functional fields has been finally verified in the paper \cite{WZ2} of Zagier and myself,  as a direct consequence of Theorem 7.2 of \cite{MR} and Theorem 2 of \cite{WZ2}. 
 \begin{thm}[Special Uniformity of Zetas; Theorem 1 of \cite{WZ2}]\label{thm1.6} For an integral  regular projective curve $X$ of genus $g$
 over $\F_q$, we have
 $$
 \ba
 \wh\zeta_{X,\F_q;n}=&\wh\zeta_{X,\F_q}^{\SL_n}(s)=q^{\binom{n}{2}(g-1)}\sum_{a=1}^{n}\sum_{\substack{k_1,\ldots,k_p>0\\ k_1+\ldots+k_p=n-a}}\frac{\wh v_{k_1}\ldots\wh v_{k_p}}{\prod_{j=1}^{p-1}(1-q^{k_j+k_{j+1})}} \frac{1}{(1-q^{ns-n+a+k_{p}})}\\
 &\hskip 2.0cm\times\wh\zeta_{X,\F_q}(ns-n+a)\sum_{\substack{l_1,\ldots,l_r>0\\ l_1+\ldots+l_r=a-1}}
 \frac{1}{(1-q^{-ns+n-a+1+l_{1}})}\frac{\wh v_{l_1}\ldots\wh v_{l_r}}{\prod_{j=1}^{r-1}(1-q^{l_j+l_{j+1})}}\\
 \ea
 $$
 where $\dis{\wh\nu_n:=\prod_{k=1}^n\wh\zeta_{X/\F_q}(k)}$ with $\,\wh\zeta_{X/\F_q}(1):=\Res_{T=1}\wh Z_{X/\F_q}(T)$.
 \end{thm}
 Indeed, in \cite{MR}, based on the theories of Hall algebra and wall-crossing, Mozgovoy-Reineke are able to  obtain a close formula for  $\wh\zeta_{X,\F_q;n}(s)$ in terms of partitions of $n$ and abelian zeta function $\wh\zeta_{X/\F_q}(s)$ of $X/\F_q$.  On the other hand, by examining the Lie structures involved in great details  in \cite{WZ2}, Zagier and myself are able to obtain the explicit formula for $\wh\zeta_{X,\F_q}^{\SL_n}(s)$ as stated in the theorem above. It is not difficult to verify that this formula of  $\wh\zeta_{X,\F_q}^{\SL_n}(s)$ coincides with the one for $\wh\zeta_{X,\F_q;n}(s)$ in  \cite{MR}. Consequently, the special uniformity of zetas for curves over finite fields is  established  successfully.

\subsection{General counting miracle}

As the first application of the special uniformity of zetas of curves over finite fields, we now deduce explicit closed formulas for the non-abelian geo-arithmetic invariants $\a_{X,\F_q;n}(mn)$ and $\b_{X,\F_q;n}(mn)$ of the curve $X$ over $\F_q$ associated to rank $n$ semi-stable vector bundles, in terms of $n,q,g$ and abelian invariants $\a_{X/\F_q}(d)$\ $(d=0,\ldots g-1)$, $\b_{X/\F_q}(0)$ and special abelian zeta values $\wh\zeta_{X/\F_q}(k)$\ $(k=1,\ldots n)$.
  Indeed, by Theorem\,\ref{thm1.2} and Theorem\,\ref{thm1.6}, we have
  $$
 \ba
 &\sum_{m=0}^{g-2}\a_{X,\F_q;n}(mn)\left(T^{m-(g-1)}+Q^{(g-1)-m}T^{(g-1)-m}\right) +\a_{X,\F_q;n}\bl n(g-1)\br+\frac{(Q-1)\b_{X,\F_q;n}(0)\cdot T}{(1-T)(1-QT)}\\
 &=q^{\binom{n}{2}(g-1)}\sum_{a=1}^{n}\sum_{\substack{k_1,\ldots,k_p>0\\ k_1+\ldots+k_p=n-a}}\frac{\,\wh v_{k_1}\ldots\wh v_{k_p}}{\prod_{j=1}^{p-1}(1-q^{k_j+k_{j+1}})} \frac{T}{(T-q^{-n+a+k_{p}})}\times\\
 &\times\left(\sum_{m=0}^{g-2}\a_{X/\F_q}(m)\left(q^{(n-a)(m-(g-1))}T^{m-(g-1)}+q^{(n-a+1)((g-1)-m)}T^{(g-1)-m}\right)+\a_{X/\F_q}\bl (g-1)\br\right.\\
&\qquad\left. +\frac{(q-1)\b_{X/\F_q}(0)\cdot q^{n-a}T}{(1-q^{n-a}T)(1-q^{n-a+1}T)}\right)\times\sum_{\substack{l_1,\ldots,l_r>0\\ l_1+\ldots+l_r=a-1}}
 \frac{1}{(1-q^{n-a+1+l_{1}}T)}\frac{\wh v_{l_1}\ldots\wh v_{l_r}}{\prod_{j=1}^{r-1}(1-q^{l_j+l_{j+1}})}.\\
\ea
 $$
 since 
 $$\ba
 \wh \zeta_{X,\F_q}(ns-n+a)=&\sum_{m=0}^{g-2}\a_{X/\F_q}(m)\left(q^{(n-a)(m-(g-1))}T^{m-(g-1)}+q^{(g-1)-m}q^{(n-a)((g-1)-m)}T^{(g-1)-m}\right)\\
&\hskip 3.0cm +\a_{X/\F_q}\bl (g-1)\br+\frac{(q-1)\b_{X/\F_q}(0)\cdot q^{n-a}T}{(1-q^{n-a}T)(1-qq^{n-a}T)}.
 \ea$$
 
To simplify our notation, set, for  each $a=1,\ldots, n$,  
 $$
\ba 
q^{-\binom{n}{2}(g-1)}\wh\zeta^{[(a)]}(s):=&\wh Z^{[(a)]}(T)=
 q^{\binom{n}{2}(g-1)}\sum_{\substack{k_1,\ldots,k_p>0\\ k_1+\ldots+k_p=n-a}}\frac{\wh v_{k_1}\ldots\wh v_{k_p}}{\prod_{j=1}^{p-1}(1-q^{k_j+k_{j+1}})} \frac{T}{(T-q^{-n+a+k_{p}})}\\
 &\times\left(\sum_{m=0}^{g-2}\a_{X/\F_q}(m)\left(q^{(n-a)(m-(g-1))}T^{m-(g-1)}+q^{(n-a+1)((g-1)-m)}T^{(g-1)-m}\right)\right.\\
&\hskip 3.0cm\left. +\a_{X/\F_q}\bl (g-1)\br+\frac{(q-1)\b_{X/\F_q}(0)\cdot q^{n-a}T}{(1-q^{n-a}T)(1-q^{n-a+1}T)}\right)\\
&\times\sum_{\substack{l_1,\ldots,l_r>0\\ l_1+\ldots+l_r=a-1}}
 \frac{1}{(1-q^{n-a+1+l_{1}}T)}\frac{\wh v_{l_1}\ldots\wh v_{l_r}}{\prod_{j=1}^{r-1}(1-q^{l_j+l_{j+1}})}.\\
 \ea
 $$
For $m=0,1,\ldots, g-1$, using the expansion $\frac{1}{1-z}=\sum_{k=0}^\infty z^k$, we have
$$\ba
&q^{-\binom{n}{2}(g-1)} \a_{X,\F_q;n}(mn)=\Res_{T=0}T^{g-2-m}\Big(\wh Z_{X,\F_q;n}(T)\Big)=\sum_{a=1}^n\Res_{T=0}\Big(T^{g-2-m}\wh Z^{[(a)]}(T)\Big)\\
% =&\sum_{a=1}^n\Res_{T=0}T^{g-2-m}\Biggl(\sum_{\substack{k_1,\ldots,k_p>0\\ k_1+
%\ldots+k_p=n-a}}\frac{\wh v_{k_1}\ldots\wh v_{k_p}}{\prod_{j=1}^{p-1}(1-q^{k_j+k_{j+1}})} (-1)\sum_{\ell=1}^\infty\bl q^{n-a-k_{p}}T\br^\ell\\
% &\times\left(\sum_{k=0}^{g-2}\a_{X/\F_q}(k)\left(q^{(n-a)(k-(g-1))}T^{k-(g-1)}+q^{(n-a+1)((g-1)-k)}T^{(g-1)-k}\right)\right.\\
%&\left. +\a_{X/\F_q}\bl (g-1)\br+(q-1)\b_{X/\F_q}(0)\sum_{\ell=1}^\infty\bl
%q^{n-a}T\br^\ell\sum_{\kappa=0}^\infty \bl q^{n-a+1}T\br^\kappa\right)\\
%&\times\sum_{\substack{l_1,\ldots,l_r>0\\ l_1+\ldots+l_r=a-1}}
 %\frac{\wh v_{l_1}\ldots\wh v_{l_r}}{\prod_{j=1}^{r-1}(1-q^{l_j+l_{j+1}})}.\sum_{\kappa=0}^\infty\bl %q^{n-a+1+l_{1}}T\br^\kappa\Biggr)\\
%=&\sum_{a=1}^n\Res_{T=0}T^{g-2-m}\Biggl(\sum_{\substack{k_1,\ldots,k_p>0\\ k_1+\ldots+k_p=n-a}}\frac{\wh v_{k_1}\ldots\wh v_{k_p}}{\prod_{j=1}^{p-1}(1-q^{k_j+k_{j+1}})} (-1)\sum_{\ell=1}^\infty\bl q^{n-a-k_{p}}T\br^\ell\\
% &\times\left(\sum_{k=0}^{g-2}\a_{X/\F_q}(k)\left(q^{(n-a)(k-(g-1))}T^{k-(g-1)}\right)+\a_{X/\F_q}\bl (g-1)\br\right)\\
%&\times\sum_{\substack{l_1,\ldots,l_r>0\\ l_1+\ldots+l_r=a-1}}
% \frac{\wh v_{l_1}\ldots\wh v_{l_r}}{\prod_{j=1}^{r-1}(1-q^{l_j+l_{j+1}})}.\sum_{\kappa=0}^\infty\bl q^{n-a+1+l_{1}}T\br^\kappa\Biggr)\\
%&\hskip 5.0cm({\rm since}\ g-2-m\geq -1)\\\ea$$ This implies that $$\ba
%q^{-\binom{n}{2}(g-1)} \a_{X,\F_q;n}(mn)
=&\sum_{a=1}^nq^{(n-a)(m-(g-1))}\Biggl(\sum_{m=k+\ell+\kappa}\sum_{k=0}^{g-2}\a_{X/\F_q}(k)\sum_{\substack{k_1,\ldots,k_p>0\\ k_1+\ldots+k_p=n-a}}\frac{\wh v_{k_1}\ldots\wh v_{k_p}}{\prod_{j=1}^{p-1}(1-q^{k_j+k_{j+1}})} (-1)\sum_{\ell=1}^\infty\bl q^{-k_{p}}\br^\ell\\
&\hskip 3.0cm\times\sum_{\substack{l_1,\ldots,l_r>0\\ l_1+\ldots+l_r=a-1}}
 \frac{\wh v_{l_1}\ldots\wh v_{l_r}}{\prod_{j=1}^{r-1}(1-q^{l_j+l_{j+1}})}.\sum_{\kappa=0}^\infty\bl q^{1+l_{1}}\br^\kappa\Biggr)+\\
&+\delta_{m,g-1}\Biggl(\a_{X/\F_q}\bl (g-1)\br\cdot\sum_{\substack{l_1,\ldots,l_r>0\\ l_1+\ldots+l_r=n-1}}
 \frac{\wh v_{l_1}\ldots\wh v_{l_r}}{\prod_{j=1}^{r-1}(1-q^{l_j+l_{j+1}})}\Biggr)\\
%&\hskip 5.0cm({\rm since}\ g-2-m\geq -1)\\
 \ea$$
 For example, if $m=0$, we get
 $$\ba
 q^{-\binom{n}{2}(g-1)}\a_{X,\F_q;n}(0)=&\bc\a_{X/\F_q}(0)\sum_{\substack{l_1,\ldots,l_r>0\\ l_1+\ldots+l_r=n-1}}
 \frac{\wh v_{l_1}\ldots\wh v_{l_r}}{\prod_{j=1}^{r-1}(1-q^{l_j+l_{j+1}})}&\qquad g\geq 2\\
\delta_{0,g-1}\Biggl(\a_{X/\F_q}\bl (g-1)\br\cdot\sum_{\substack{l_1,\ldots,l_r>0\\ l_1+\ldots+l_r=n-1}}
 \frac{\wh v_{l_1}\ldots\wh v_{l_r}}{\prod_{j=1}^{r-1}(1-q^{l_j+l_{j+1}})}\Biggr)&\qquad g=1\\
\ec \ea$$
That is to say, we have proved the following counting miracle relation, since 
\be\a_{X/\F_q}(0)=\sum_{L\in \Pic^0(X)}\frac{q^{h^0(X,L)}-1}{q-1}=\frac{q^{h^0(X,\O_X)}-1}{q-1}=1.\ee

\begin{thm}[Counting Miracle]
For an integral regular projective curve $X$ of genus $g\geq 1$ over a finite field $\F_q$, 
\be
q^{-\binom{n}{2}(g-1)} \a_{X,\F_q;n}(0)=\sum_{\substack{l_1,\ldots,l_r>0\\ l_1+\ldots+l_r=n-1}}
 \frac{\wh v_{l_1}\ldots\wh v_{l_r}}{\prod_{j=1}^{r-1}(1-q^{l_j+l_{j+1}})}=q^{-\binom{n-1}{2}(g-1)}\b_{X,\F_q,n-1}(0)
\ee
\end{thm}

The counting miracle was first conjectured in \cite{W}. It is established in \cite{WZ1} for elliptic curves with a heavy use of combinatorial technique in 2014, after examining Atiyah bundles in details. In September 2016, adopting a totally different method, K. Sugahara established this counting miracle. It was later reverified in \cite{MR} independently using Hall algebra and wall crossing.  In an appendix to this subsection, Sugahara's proof will be presented.

Similarly,
$$
 \ba
&q^{-\binom{n}{2}(g-1)} \Big(\b_{X,\F_q;n}(0)+Q\a_{X,\F_q;n}\bl(g-2)n\br\Big)=
\Res_{T=0}T^{-2}\wh Z_{X,\F_q;n}(T)\\
%=&\sum_{a=1}^n\Res_{T=0}T^{-2}\Biggl(\sum_{\substack{k_1,\ldots,k_p>0\\ k_1+\ldots+k_p=n-a}}\frac{\wh v_{k_1}\ldots\wh v_{k_p}}{\prod_{j=1}^{p-1}(1-q^{k_j+k_{j+1}})} (-1)\sum_{\ell=1}^\infty\bl q^{n-a-k_{p}}T\br^\ell\\
 %&\times\left(\sum_{k=0}^{g-2}\a_{X/\F_q}(k)\left(q^{(n-a)(k-(g-1))}T^{k-(g-1)}+q^{(n-a+1)((g-1)-k)}T^{(g-1)-k}\right)\right.\\
%&\left. +\a_{X/\F_q}\bl (g-1)\br+(q-1)\b_{X/\F_q}(0)\sum_{\ell=1}^\infty\bl
%q^{n-a}T\br^\ell\sum_{\kappa=0}^\infty \bl q^{n-a+1}T\br^\kappa\right)\\
%&\times\sum_{\substack{l_1,\ldots,l_r>0\\ l_1+\ldots+l_r=a-1}}
 %\frac{\wh v_{l_1}\ldots\wh v_{l_r}}{\prod_{j=1}^{r-1}(1-q^{l_j+l_{j+1}})}.\sum_{\kappa=0}^\infty\bl q^{n-a+1+l_{1}}T\br^\kappa\Biggr)\\
% =&\sum_{a=1}^n\Res_{T=0}T^{-2}\Biggl(\sum_{\substack{k_1,\ldots,k_p>0\\ k_1+\ldots+k_p=n-a}}\frac{\wh v_{k_1}\ldots\wh v_{k_p}}{\prod_{j=1}^{p-1}(1-q^{k_j+k_{j+1}})} (-1)\sum_{\ell=1}^\infty\bl q^{n-a-k_{p}}T\br^\ell\\
 %&\times\left(\sum_{k=0}^{g-2}\a_{X/\F_q}(k)\left(q^{(n-a)(k-(g-1))}T^{k-(g-1)}\right)+\a_{X/\F_q}\bl (g-1)\br+(q-1)\b_{X/\F_q}(0)\sum_{\ell=1}\bl
%q^{n-a}T\br^\ell\right)\\
%&\times\sum_{\substack{l_1,\ldots,l_r>0\\ l_1+\ldots+l_r=a-1}}
 %\frac{\wh v_{l_1}\ldots\wh v_{l_r}}{\prod_{j=1}^{r-1}(1-q^{l_j+l_{j+1}})}.\sum_{\kappa=0}^\infty\bl q^{n-a+1+l_{1}}T\br^\kappa\Biggr)\\
 =&\sum_{a=1}^nq^{n-a}\Biggl(-\sum_{\substack{k_1,\ldots,k_p>0\\ k_1+\ldots+k_p=n-a}}\frac{\wh v_{k_1}\ldots\wh v_{k_p}}{\prod_{j=1}^{p-1}(1-q^{k_j+k_{j+1}})} \sum_{\substack{l_1,\ldots,l_r>0\\ l_1+\ldots+l_r=a-1}}
 \frac{\wh v_{l_1}\ldots\wh v_{l_r}}{\prod_{j=1}^{r-1}(1-q^{l_j+l_{j+1}})}\\
 &\hskip 4.0cm\times \sum_{k+\ell+\kappa=g}\sum_{k=0}^{g-2}\a_{X/\F_q}(k)\sum_{\ell=1}^g\bl q^{-k_{p}}\br^\ell\sum_{\kappa=0}^g\bl q^{1+l_{1}}\br^\kappa\\
 &\hskip 1.40cm-\a_{X/\F_q}\bl (g-1)\br\cdot\sum_{\substack{k_1,\ldots,k_p>0\\ k_1+\ldots+k_p=n-a}}\frac{\wh v_{k_1}\ldots\wh v_{k_p}}{\prod_{j=1}^{p-1}(1-q^{k_j+k_{j+1}})} \bl q^{-k_{p}}\br\\
 &\hskip 6.0cm\times \sum_{\substack{l_1,\ldots,l_r>0\\ l_1+\ldots+l_r=a-1}}
 \frac{\wh v_{l_1}\ldots\wh v_{l_r}}{\prod_{j=1}^{r-1}(1-q^{l_j+l_{j+1}})}\\
 &\hskip 1.40cm+(q-1)\b_{X/\F_q}(0)\cdot\sum_{\substack{l_1,\ldots,l_r>0\\ l_1+\ldots+l_r=n-1}}
 \frac{\wh v_{l_1}\ldots\wh v_{l_r}}{\prod_{j=1}^{r-1}(1-q^{l_j+l_{j+1}})}
 \Biggr)\\
\ea
 $$
 This then gives a closed formula for $\b_{X,\F_q;n}(0)$ after substracting $Q\cdot \a_{X,\F_q;n}\bl(g-2)n\br$ obtained above. Thus all in all, we have proved the following

 \begin{thm}[General Counting Miracle]\label{thm1.7} For an integral regular projective curve $X$ of genus $g$ on $\F_q$, its non-abelian invariants $\a_{X,\F_q,n}(mn)\ (0\leq  m\leq g-1)$ and $\b_{X,\F_q;n}(0)$ for semi-stable vector bundles of rank $n$ are given by
 $$\ba
&q^{-\binom{n}{2}(g-1)} \a_{X,\F_q;n}(mn)\\
=&\sum_{a=1}^nq^{(n-a)(m-(g-1))}\Biggl(
-\sum_{\substack{k_1,\ldots,k_p>0\\ k_1+\ldots+k_p=n-a}}\frac{\wh v_{k_1}\ldots\wh v_{k_p}}{\prod_{j=1}^{p-1}(1-q^{k_j+k_{j+1}})} \sum_{\substack{l_1,\ldots,l_r>0\\ l_1+\ldots+l_r=a-1}}
 \frac{\wh v_{l_1}\ldots\wh v_{l_r}}{\prod_{j=1}^{r-1}(1-q^{l_j+l_{j+1}})}\\
 &\hskip 5.0cm\times\sum_{m=k+\ell+\kappa}\sum_{k=0}^{g-2}\a_{X/\F_q}(k)\sum_{\ell=1}^\infty\bl q^{-k_{p}}\br^\ell\sum_{\kappa=0}^\infty\bl q^{1+l_{1}}\br^\kappa\Biggr)\\
&\quad+\delta_{m,g-1}\Biggl(\a_{X/\F_q}\bl (g-1)\br\cdot\sum_{\substack{l_1,\ldots,l_r>0\\ l_1+\ldots+l_r=n-1}}
 \frac{\wh v_{l_1}\ldots\wh v_{l_r}}{\prod_{j=1}^{r-1}(1-q^{l_j+l_{j+1}})}\Biggr)\\
 \ea$$
 and
 $$\ba
&q^{-\binom{n}{2}(g-1)} \Big(\b_{X,\F_q;n}(0)\Big)\\
 =&\sum_{a=1}^n\Biggl(-\sum_{\substack{k_1,\ldots,k_p>0\\ k_1+\ldots+k_p=n-a}}\frac{\wh v_{k_1}\ldots\wh v_{k_p}}{\prod_{j=1}^{p-1}(1-q^{k_j+k_{j+1}})} \sum_{\substack{l_1,\ldots,l_r>0\\ l_1+\ldots+l_r=a-1}}
 \frac{\wh v_{l_1}\ldots\wh v_{l_r}}{\prod_{j=1}^{r-1}(1-q^{l_j+l_{j+1}})}\\
 &\hskip 3.0cm\times \left(q^{n-a}\sum_{k+\ell+\kappa=g}-q^a\sum_{k+\ell+\kappa=g-2}\right)\left(\sum_{k=0}^{g-2}\a_{X/\F_q}(k)\sum_{\ell=1}^\infty\bl q^{-k_{p}}\br^\ell\sum_{\kappa=0}^\infty\bl q^{1+l_{1}}\br^\kappa\right)\\
 &\quad-\a_{X/\F_q}\bl (g-1)\br\cdot\sum_{\substack{k_1,\ldots,k_p>0\\ k_1+\ldots+k_p=n-a}}\frac{\wh v_{k_1}\ldots\wh v_{k_p}}{\prod_{j=1}^{p-1}(1-q^{k_j+k_{j+1}})} \bl q^{n-a-k_{p}}\br \sum_{\substack{l_1,\ldots,l_r>0\\ l_1+\ldots+l_r=a-1}}
 \frac{\wh v_{l_1}\ldots\wh v_{l_r}}{\prod_{j=1}^{r-1}(1-q^{l_j+l_{j+1}})}\\
 &\hskip 5.0cm+(q-1)q^{n-a}\b_{X/\F_q}(0)\cdot\sum_{\substack{l_1,\ldots,l_r>0\\ l_1+\ldots+l_r=n-1}}
 \frac{\wh v_{l_1}\ldots\wh v_{l_r}}{\prod_{j=1}^{r-1}(1-q^{l_j+l_{j+1}})}
 \Biggr)\\
\ea
 $$
 \end{thm}
Recall that
\be
\b_{X,\F_q;n}(mn)=\b_{X,\F_q;n}(0)\quad(\forall m\in \Z)\qan \a_{X,\F_q;n}(mn)=q^{m(n-(g-1)}\b_{X,\F_q;n}(0)\quad(m\geq g)
\ee
this theorem in fact gives all the values of $\a_{X,\F_q;n}(mn)$ and $\b_{X,\F_q;n}(mn)$ for all $m\in \Z$, since easily
\be
\a_{X,\F_q;n}(mn)=0\qquad(m<0).
\ee
We point out in passing when $n$ does not divide $d$,\footnote{In \cite{MR}, Mozegovoy-Reineke call such pairs $(n,d)$ generic  and provide a method to calculated the $\a$ and $\b$ invariants when $(n,d)$ are generic.} the value of $\a_{X,\F_q;n}(d)$  
 and $\b_{X,\F_q;n}(d)$ have been obtained in \cite{MR} and \cite{Z}, respectively.

 We end this subsection with the following comments on $\wh Z_{X,\F_q}^{\SL_n}(T)$. By the special uniformity, this function is equal to the rank $n$ zeta function of $X$, which itself is a rational function of the form
 \be
 \wh Z_{X,\F_q;n}(T)=\frac{P_{X,\F_q;n}(T)}{T^{g-1}(1-T)(1-QT)}.
 \ee
Here $P_{X,\F_q;n}(T)$ is a polynomial of degree $2g$ in $T$ with real coefficients.
However, in the summand of $\wh\zeta^{[(a)]}(s)$,
from  the first group $\frac{\wh v_{k_1}\ldots\wh v_{k_p}}{\prod_{j=1}^{p-1}(1-q^{k_j+k_{j+1}})} \frac{T}{(T-q^{-n+a+k_{p}})}$, particularly the term $\frac{T}{(T-q^{-n+a+k_{p}})}=\frac{q^{(k_1+\ldots+k_{p-1})}T}{q^{(k_1+\ldots+k_{p-1})}T-1}$, we see that the denominators are given by 
\be\frac{1}{T-1},\ \frac{1}{qT-1},\ \ldots,\ \frac{1}{q^{n-a-1}T-1}\ee
In parallel, from the second of the whole bracket, particularly the term $\frac{(q-1)\b_{X/\F_q}(0)\cdot q^{n-a}T}{(1-q^{n-a}T)(1-q^{n-a+1}T)}$. we see that the denominators are given by 
\be
\frac{1}{q^{n-a}T-1},\ \frac{1}{q^{n-a+1}T-1},
\ee
Similarly, from the third group, particularly the term $\frac{1}{(1-q^{n-a+1+l_{1}}T)}$, we see that the denominators are given by 
\be
\frac{1}{q^{n-a+1+1}T-1},\ \frac{1}{(1-q^{n-a+1+2}T)},\ \ldots,\ \frac{1}{(1-q^{n-a+1+a-1}T)}
\ee
Consequently,  
$$\ba \xi^{[(a)]}(s):=&\left(T^{g-1}\prod_{\ell=0}^n(q^\ell T-1)\right)\cdot \wh\zeta^{[(a)]}(s)\\
%= &\sum_{\substack{k_1,\ldots,k_p>0\\ k_1+\ldots+k_p=n-a}}\frac{\wh v_{k_1}\ldots\wh v_{k_p}}%{\prod_{j=1}^{p-1}(1-q^{k_j+k_{j+1}})} \frac{q^{k_1+\ldots+k_{p-1}}T}{(q^{k_1+\ldots+k_{p-1}}T-1)}%\prod_{\ell=0}^{n-a-1}(q^\ell T-1)\\
 %&\times\Biggl(\Big(\sum_{m=0}^{g-2}\a_{X/\F_q}(m)\left(q^{(n-a)(m-(g-1))}T^{m}+q^{(n-a+1)((g-1)-m)}T^{2(g-1)-m}\right)\\
 %&\hskip 2.0cm+\a_{X/\F_q}\bl (g-1)\br T^{g-1}\Big){(1-q^{n-a}T)(1-q^{n-a+1}T)}\\
%&\hskip 5.0cm+(q-1)\b_{X/\F_q}(0)\cdot q^{n-a}T^g\Biggr)\\
%&\times\sum_{\substack{l_1,\ldots,l_r>0\\ l_1+\ldots+l_r=a-1}}
 %\frac{1}{(1-q^{n-a+1+l_{1}}T)}\frac{\wh v_{l_1}\ldots\wh v_{l_r}}{\prod_{j=1}^{r-1}(1-q^{l_j+l_{j+1}})}\prod_{\ell=n-a+2}^n(q^\ell T-1)T^{g-1}.\\
% \ea$$
 %$$\ba
 = &\sum_{\substack{k_1,\ldots,k_p>0\\ k_1+\ldots+k_p=n-a}}\frac{\wh v_{k_1}\ldots\wh v_{k_p}}{\prod_{j=1}^{p-1}(1-q^{k_j+k_{j+1}})} q^{k_1+\ldots+k_{p-1}}T\prod_{\substack{0\leq\ell\leq n-a-1\\ \ell\not=k_1+\ldots+k_{p-1}}}(q^\ell T-1)\\
\qquad &\times\Biggl(\Big(\sum_{m=0}^{g-2}\a_{X/\F_q}(m)\left(q^{(n-a)(m-(g-1))}T^{m}+q^{(n-a+1)((g-1)-m)}T^{2(g-1)-m}\right)\\
 &\hskip 2.0cm+\a_{X/\F_q}\bl (g-1)\br T^{g-1}\Big){(1-q^{n-a}T)(1-q^{n-a+1}T)}\\
&\hskip 5.0cm\left. +(q-1)\b_{X/\F_q}(0)\cdot q^{n-a}T^g\right.\Biggr)\\
\qquad\qquad&\times\sum_{\substack{l_1,\ldots,l_r>0\\ l_1+\ldots+l_r=a-1}}
\frac{\wh v_{l_1}\ldots\wh v_{l_r}}{\prod_{j=1}^{r-1}(1-q^{l_j+l_{j+1}})}\prod_{\substack{n-a+2\leq \ell\leq n\\\ell\not=n-a+1+l_1}}(q^\ell T-1).\\
\ea
$$
becomes a polynomial of degree $(n-a)+2g+(a-2)=n+2(g-1)$ which is independent of $a$. Therefore, 
$\wh Z_{X,\F_q}^{\SL_n}(T)$ is a rational function of the form
\be
\frac{P_{X,\F_q}^{\SL_n}(T)}{T^{g-1}\prod_{\ell=0}^n(q^\ell T-1)}
\ee
where $P_{X,\F_q}^{\SL_n}(T)$ is a polynomial of degree $n+2(g-1)$ in $T$.
By comparing this with $\frac{P_{X,\F_q;n}(T)}{T^{g-1}\prod_{\ell=0}^n(q^\ell T-1)}$, we see that in fact there are massive cancellations among the $Z^{[a]}(s)$ when combining the summation $\sum_{a=1}^n Z^{[a]}(T)$ to obtain $\wh Z_{X,\F_q}^{\SL_n}(T)$ such that among the product $\prod_{\ell=0}^{n}(q^\ell T-1)$, all the factors $(q T-1),\ (q^2T-1),\ldots, (q^{n-1}T-1)$ will be finally cancelled out from the numerator (so as to leave only the factor $T^{g-1}(1-T)(1-QT)$ in the denominator). This is one of the reasons why
the Riemann Hypothesis for high rank zeta functions of curves over finite fields becomes quite complicated, even comparing with what has happened for high rank zeta functions of  number fields.

\subsection*{Appendix:}
\centerline {\bf\Large Counting Miracle} {~}\\
\centerline{\bf Kotaro Sugahara}\footnote{When this note was completed in September, 2016, I was a PhD student in Graduate School of Mathematics, Kyushu university.}
{~}\\
\centerline{\bf   Fuji Soft}{~}\\

\addcontentsline{toc}{subsection}{Appendix A}

Let $X$ be an integral regular projective curve of genus $g$ over $\F_q$, the finite field with $q$ elements. 
For any coherent sheaves $\mathcal{A}$, $\mathcal{B}$ and $\mathcal{E}$ on $X$, for our own use, we introduce the auxiliary spaces
$$\begin{aligned}
\rm{Fil}(\mathcal{A}, \mathcal{B}; \mathcal{E}):=&\Big\{0\subset \mathcal{F} \subset \mathcal{E}: \mathcal{F}\simeq \mathcal{A},\  \mathcal{E}/\mathcal{F}\simeq \mathcal{B}\Big\},\\
\rm{W}(\mathcal{A}, \mathcal{B}; \mathcal{E}):=&\Big\{(\textit{f, g})\in \rm{Hom}(\mathcal{A}, \mathcal{E})\times \rm{Hom}(\mathcal{E}, \mathcal{B}): \textit{f\ injective, g\ surjective, g}\circ \textit{f}=0\Big\},\\
\rm{U}(\mathcal{A}, \mathcal{B}; \mathcal{E}):=&\Big\{\textit{f}\in \rm{Hom}(\mathcal{A}, \mathcal{E}): \textit{f\ injective, }\exists\,{\rm surjective}\, \textit {g}\in \rm{Hom}(\mathcal{E}, \mathcal{B})\textit{ s.t. g}\circ \textit{f}=0\Big\},\\\end{aligned}$$ 
the associated morphisms
$$\begin{aligned}
\varphi:\rm{W}(\mathcal{A}, \mathcal{B}; \mathcal{E})&\rightarrow \rm{Fil}(\mathcal{A}, \mathcal{B}; \mathcal{E});\qquad (\textit{f, g})\mapsto (0\subset {\rm{Im}}\, ( f )\subset \mathcal{E})\\
\psi:\rm{W}(\mathcal{A}, \mathcal{B}; \mathcal{E})&\rightarrow \rm{Ext}^{1}(\mathcal{B}, \mathcal{A});\qquad (\textit{f, g})\mapsto (0\rightarrow \mathcal{A}\xrightarrow{\textit{f}} \mathcal{E}\xrightarrow{\textit{g}} \mathcal{B}\rightarrow 0)\\
\pi:\rm{W}(\mathcal{A}, \mathcal{B}; \mathcal{E})&\rightarrow \rm{U}(\mathcal{A}, \mathcal{B}; \mathcal{E});\qquad (\textit{f, g})\mapsto \textit{f}\\\end{aligned}$$ 
and the following natural actions:
$$\begin{aligned}\chi:\{(\rm{Aut}\,\mathcal{A})^{op}&\times \rm{Aut}\,\mathcal{B}\} \times \rm{W}(\mathcal{A}, \mathcal{B}; \mathcal{E})\rightarrow \rm{W}(\mathcal{A}, \mathcal{B}; \mathcal{E});\ \  (\rho, \sigma,(\textit{f, g}))\mapsto (\textit{f}\circ \rho, \sigma \circ \textit{g})\\
\mu:\rm{Aut}\,\mathcal{E}\times \rm{W}&(\mathcal{A}, \mathcal{B}; \mathcal{E})\rightarrow \rm{W}(\mathcal{A}, \mathcal{B}; \mathcal{E});\qquad (\tau, (\textit{f, g}))\mapsto (\tau\circ \textit{f, g} \circ \tau^{-1})
\end{aligned}$$  

Accordingly, we have \\

\noindent
(i) $\varphi$ is surjective, and 

\noindent
(ii) There is a natural bijection between the fiber of $\varphi$ and $(\rm{Aut}\,\mathcal{A})^{op}\times \rm{Aut}\,\mathcal{B}$. \\

This is a direct consequence of the five lemma. Indeed, this follows from the facts that $\varphi(f, g)=\varphi(f', g')$ if and only if there exists $(\rho, \sigma)\in (\rm{Aut}\,\mathcal{A})^{op}\times \rm{Aut}\,\mathcal{B}$ such that $(f', g')=\chi(\rho, \sigma, (f, g))$, and that the action $\chi$ is free. Hence, we have
\begin{eqnarray}\label{eq31}
\# \rm{W}(\mathcal{A}, \mathcal{B}; \mathcal{E})\big/\big(\# (\rm{Aut}\,\mathcal{A})\cdot \# \rm{Aut}\,\mathcal{B}\big)=\# \rm{Fil}(\mathcal{A}, \mathcal{B}; \mathcal{E}).
\end{eqnarray}   

\noindent
(iii) The image of $\psi$ is exactly the set $\rm{Ext}_{\mathcal{E}}^{1}(\mathcal{B}, \mathcal{A})$ of isomorphism classes of extensions of $\mathcal{B}$ by $\mathcal{A}$  the middle term of which is isomorphic to $\mathcal{E}$.
Moreover,  any fiber of $\psi$ is an orbit of $\rm{W}(\mathcal{A}, \mathcal{B}; \mathcal{E})$ under the action $\mu$. 

Indeed, this follows from the fact that $\psi(f,g)=\psi(f',g')$ if and only if there exists $\tau\in \rm{Aut}\,\mathcal{E}$ such that $(f', g')=\mu(\tau,  (f, g))$. In addition,\\

\noindent
(iv) The stabilizer group of $(f, g)\in \rm{W}(\mathcal{A}, \mathcal{B}; \mathcal{E})$ under the action $\mu$ is isomorphic to $\rm{Hom}(\mathcal{B},\mathcal{A}).$

Indeed, for any element $\tau\in\Aut\,\cE$, if we write it as $\tau=\begin{pmatrix}\tau_{\aA}&\tau_{\aA\bB}\\
\tau_{\bB\aA}&\tau_{\bB}\end{pmatrix}$ with
$$\tau_{\aA}\in \Aut\,\aA,\ \tau_{\bB}\in \Aut\,\bB,\ \tau_{\aA\bB}\in \Hom(\aA,\bB),\ 
\tau_{\bB\aA}\in \Hom(\bB,\aA),$$
then a direct calculation shows that $f=\tau\circ f$ if and only if 
$\tau=\begin{pmatrix}\mathrm{Id}_{\aA}&0\\
\tau_{\bB\aA}&\tau_{\bB}\end{pmatrix}.$
Here we have used  the claim that $\rm{Hom}(\mathcal{A}, \mathcal{B})$ part $\tau_{\aA\bB}$ must be 0 since we have the inclusion
$\tau(f(\mathcal{A}))\subset f(\mathcal{A})$ by condition. Similarly, 
$g=g\circ \tau^{-1}$ if and only if $g\circ\tau=g$ if and only if
$\tau=\begin{pmatrix}\tau_{\aA}&0\\
\tau_{\bB\aA}&\mathrm{Id}_{\bB}\end{pmatrix}.$
Thus, $(f, g)=\mu(\tau,  (f, g))$ if and only if 

\begin{equation*}
\tau =
\begin{pmatrix}
\rm{Id}_{\mathcal{A}} & 0\\
\tau_{\mathcal{B} \mathcal{A}} & \rm{Id}_{\mathcal{B}}
\end{pmatrix}\in \begin{pmatrix}
\rm{Id}_{\mathcal{A}} & 0\\
\Hom({\mathcal{B}, \mathcal{A}}) & \rm{Id}_{\mathcal{B}}
\end{pmatrix}
.
\end{equation*}
 Consequently, the stabilizer group of $(f, g)\in \rm{W}(\mathcal{A}, \mathcal{B}; \mathcal{E})$ under the action $\mu$ is isomorphic to $\rm{Hom}(\mathcal{B},\mathcal{A}).$

Therefore, we obtain the following relation
\begin{eqnarray}\label{eq32}
\# \rm{W}(\mathcal{A}, \mathcal{B}; \mathcal{E})\big/\big(\# \rm{Aut}\mathcal{E}\big/\# (\rm{Hom}(\mathcal{B}, \mathcal{A}))\big)=\# \rm{Ext}_{\mathcal{E}}^{1}(\mathcal{B}, \mathcal{A}).
\end{eqnarray}   

By \eqref{eq31} and \eqref{eq32}, we have
\begin{eqnarray}\label{eq33}
\# \rm{Fil}(\mathcal{A}, \mathcal{B}; \mathcal{E})=\frac{\# \rm{Ext}_{\mathcal{E}}^{1}(\mathcal{B}, \mathcal{A})\cdot \# \rm{Aut}\mathcal{E}}{\# \rm{Aut}\,\mathcal{A}\cdot \# \rm{Aut}\mathcal{B}\cdot \# (\rm{Hom}(\mathcal{B}, \mathcal{A}))}.
\end{eqnarray}   

We mention in passing that this relation is an analogue of the formula given by [Ri]. \\

\noindent
(v) $\pi$ is surjective, and each fiber of $\pi$ is isomorphic to $\rm{Aut}\,\mathcal{B}$. 

This follows from the fact that $\pi(f, g)=\pi(f',g')$ if and only if $f=f'$ and there exists $\sigma \in \rm{Aut}\,\mathcal{B}$, by the five lemma. Consequently, we have
\begin{eqnarray}\label{eq34}
\# \rm{W}(\mathcal{A}, \mathcal{B}; \mathcal{E})\big/\# \rm{Aut}\,\mathcal{B}=\# \rm{U}(\mathcal{A}, \mathcal{B}; \mathcal{E}).
\end{eqnarray}   

\noindent
{\bf Thorem.}(Counting Miracle) {\it
Let $\mathcal{E}_{0}$ be a stable vector bundle of rank $m\ (< n)$ and degree $0$.
Then}
\begin{eqnarray*}
\sum_{\mathcal{E}\in \mathcal{M}_{X,n}(0)}\frac{q^{\#\Hom(\mathcal{E}_{0}, \mathcal{E})}-1}{\# \rm{Aut}\, \mathcal{E}}=q^{m(n-m)(g-1)}\sum_{\mathcal{F}\in \mathcal{M}_{X,n-m}(0)}\frac{1}{\# \rm{Aut}\, \mathcal{F}}.
\end{eqnarray*}
{\it In particular,}
\begin{eqnarray*}
\alpha_{X,\F_q;n+1}(0)=q^{n(g-1)}\beta_{X,\F_q;n}(0).
\end{eqnarray*}

This counting miracle was  conjectured by Weng in \cite{W}. It was first established  by Weng and Zagier \cite{WZ1} for elliptic curves,  based on basic properties of Atiyah bundles and heavy combinatorial techniques. Our method here is totally different.

\bp It suffices to prove the first relation, since the second comes directly by applying the first to $\cE_0=\O_X$ and $m=1$. 
To prove the first, by the stability of $\cE_0$
\begin{eqnarray*}
\sum_{\mathcal{E}\in \mathcal{M}_{X,n}(0)}\frac{q^{\#\Hom(\mathcal{E}_{0}, \mathcal{E})}-1}{\# \rm{Aut}\, \mathcal{E}}
&=&\sum_{\mathcal{E}\in \mathcal{M}_{X,n}(0)}\frac{\# \{f\in \Hom(\mathcal{E}_{0}, \mathcal{E}) : f \ {\textrm {injective}}\}}{\# \rm{Aut}\, \mathcal{E}}\\
%&&\hskip 3.0cm (\textrm{by the stability of}\ \mathcal{E}_{0})\\
&=&\sum_{\mathcal{E}\in \mathcal{M}_{X,n}(0)}\frac{1}{\# \rm{Aut}\, \mathcal{E}}\,
\# \rm{U}(\mathcal{E}_{0}, \mathcal{E}/\mathcal{E}_0; \mathcal{E})\\
\end{eqnarray*}
since the category of semi-stable bundles of degree  0 is abelian.
Hence, by \eqref{eq34} and \eqref{eq33}, this latest quantity can be written as
\begin{eqnarray*}
&=&\sum_{\mathcal{E}\in \mathcal{M}_{X,n}(0)} \frac{1}{\# \rm{Aut}\, \mathcal{E}}\cdot \frac{\# \rm{W}(\mathcal{E}_{0}, \mathcal{E}/\mathcal{E}_0; \mathcal{E})}{\# \rm{Aut}(\mathcal{E}/\mathcal{E}_0)}\\
&=&\sum_{\mathcal{E}\in \mathcal{M}_{X,n}(0)}\frac{1}{\# \rm{Aut}\, \mathcal{E}}\cdot \frac{\# \rm{Ext}_{\mathcal{E}}^{1}(\mathcal{E}/\mathcal{E}_0, \mathcal{E}_{0})\cdot \# \rm{Aut} \mathcal{E}}{\# \rm{Aut} \,(\mathcal{E}/\mathcal{E}_0)\cdot \# \rm{Hom}(\mathcal{E}/\mathcal{E}_0, \mathcal{E}_{0})}\\
%&=&\sum_{\mathcal{F}\in \mathcal{M}_{X,n-m}(0)}\frac{\# \rm{Ext}_{\mathcal{E}}^{1}(\mathcal{F}, \mathcal{E}_{0})}{\# \rm{Aut} \,\mathcal{F}\cdot \# \rm{Hom}(\mathcal{F}, \mathcal{E}_{0})}\\
&=&\sum_{\mathcal{F}\in \mathcal{M}_{X,n-m}(0)}\frac{\# \rm{Ext}^{1}(\mathcal{F}, \mathcal{E}_{0})}{\# \rm{Aut} \,\mathcal{F}\cdot \# \rm{Hom}(\mathcal{F}, \mathcal{E}_{0})}\\
\end{eqnarray*}
by changing from $\mathcal E$ to $\mathcal F$, since the category of semistable bundles of degree  0  is  abelian.
Therefore, by the Riemann-Roch theorem
\begin{eqnarray*}
\sum_{\mathcal{E}\in \mathcal{M}_{X,n}(0)}\frac{q^{\#\Hom(\mathcal{E}_{0}, \mathcal{E})}-1}{\# \rm{Aut}\, \mathcal{E}}
&=&q^{m(n-m)(g-1)}\sum_{\mathcal{F}\in \mathcal{M}_{X,n-m}(0)}\frac{1}{\# \rm{Aut}\, \mathcal{F}}
\end{eqnarray*}
as wanted.
\ep

\section{Riemann hypothesis for rank two zeta: Yoshida's approach}

Applying Theorem\,\ref{thm1.6} to $n=2$, easily we conclude that, up to a constant fact depending only on the genus $g$ of the curve $X/\F_q$, which certainly does not affect our discussion on zeta zeros,
\be
\wh\zeta^{\,\SL_2}_{X,\F_q}(s)=\frac{\wh\zeta^{~}_{X/\F_q}(2s)}{1-q^{2-2s}}-\frac{\wh\zeta^{~}_{X/\F_q}(2s-1)}{q^{2s}-1}
\ee
We first use the functional equation to obtain
\be
\wh\zeta^{\,\SL_2}_{X,\F_q}(s)=\frac{\wh\zeta^{~}_{X/\F_q}(1-2s)}{1-q^{2-2s}}-\frac{\wh\zeta_{X,/F_q}(2s-1)}{q^{2s}-1}
=\frac{\wh\zeta^{~}_{X/\F_q}(2\s)}{1-q^{1+2\s}}-\frac{\wh\zeta^{~}_{X/\F_q}(-2\s)}{q^{1-2\s}-1}
\ee
where $s=:\frac{1}{2}-\s$.
Therefore,
\be\label{eq7}
\wh\zeta^{\,\SL_2}_{X,\F_q}(s)=0\qquad{\rm
if\ and\ only\ if}\qquad
\frac{\wh\zeta^{~}_{X/\F_q}(2\s)}{1-q^{1+2\s}}=\frac{\wh\zeta^{~}_{X/\F_q}(-2\s)}{q^{1-2\s}-1}.
\ee
Write now 
\be
\wh\zeta^{~}_{X/\F_q}(s)=\frac{\prod_{i=1}^g(1-\om_iq^{-s})(1-\ov\om_iq^{-s})}{q^{-s(g-1)}(1-q^{-s})(1-q^{1-s})}
\ee
where $\om_i\in \C\smm\R$ and $|\om_i|=\sqrt q$, guaranteed by the Hasse-Weil Theorem, or better, the Riemann hypothesis for the Artin zeta function of $X/\F_q$. In particular, 
\be
|\om_i+\ov\om_i|\leq 2\sqrt q<q+1\qquad(\forall 1\leq i\leq g).
\ee
Accordingly, \eqref{eq7} becomes
\be
q^{4\s(g-1)}(q^{1-2\s}-1)\frac{\prod_{i=1}^g(1-\om_iq^{-2\s})(1-\ov\om_iq^{-2\s})}{(1-q^{-2\s})(1-q^{1-2\s})}
=(1-q^{1+2\s})\frac{\prod_{i=1}^g(1-\om_iq^{2\s})(1-\ov\om_iq^{2\s})}{(1-q^{2\s})(1-q^{1+2\s})}
\ee
This is equivalent to
\be
q^{4\s(g-1)}(1-q^{2\s})\prod_{i=1}^g(1-\om_iq^{-2\s})(1-\ov\om_iq^{-2\s})=(1-q^{-2\s})\prod_{i=1}^g(1-\om_iq^{2\s})(1-\ov\om_iq^{2\s})
\ee
In particular, we should have
\be\label{eq12}
|q^{4\s(g-1)}|\cdot |(1-q^{2\s})|\cdot \prod_{i=1}^g|(1-\om_iq^{-2\s})(1-\ov\om_iq^{-2\s})|=|(1-q^{-2\s})|\cdot \prod_{i=1}^g|(1-\om_iq^{2\s})(1-\ov\om_iq^{2\s})|
\ee
\begin{lem}[Yoshida]\label{lem2.1} Fix a real number $q>1$. Let $\a,\,\b\in \C$ and write $c=\a+\b$. Assume that $\a\b=q$ and that $c\in \R$ satisfies $|c|\leq q+1$. Then for $w\in \C$, we have
\be
|w-\a|\cdot|w-\b|\ \bc
&\hskip -0.3cm>|1-\a w|\cdot|1-\b w|\qquad{\rm if}\quad |w|<1\\
&\hskip -0.3cm <|1-\a w|\cdot|1-\b w|\qquad{\rm if}\quad |w|>1.
\ec
\ee
\end{lem}
When $\a=1$ and $\b=q$,  this lemma degenerates to an estimate on the fractional transformation $T_q(w):=\dis{\frac{w-q}{1-qw}}$. In this sense, Yoshida's lemma is a natural degree 2 generalization of that for fractional transformations. Even an elementary  proof can be given immediately, we delay the details till the proof of Lemma\,\ref{lem3.4} below, which itself is a generalization of Yoshida's lemma.
 
Therefore, for the left hand side of \eqref{eq12}, by applying Lemma\,\ref{lem2.1}, we have
$$\ba
&|q^{4\s(g-1)}|\cdot |(1-q^{2\s})|\cdot \prod_{i=1}^g|(1-\om_iq^{-2\s})(1-\ov\om_iq^{-2\s})|\\
=&|(1-q^{2\s})||q^{-4\s}|\cdot \prod_{i=1}^g|(q^{2\s}-\om_i)(q^{2\s}-\ov\om_i)|\\
&\bc
&\hskip -0.3cm>\prod_{i=1}^g|(1-\om_i q^{2\s})(1-\ov\om_i q^{2\s})|\cdot |(q^{-2\s}-1)||q^{(1-2g)2\s}|\qquad{\rm if}\quad \Re(\s)<0\\
&\hskip -0.3cm <\prod_{i=1}^g|(1-\om_i q^{2\s})(1-\ov\om_i q^{2\s})|\cdot|(q^{-2\s}-1)||q^{(1-2g)2\s}|\qquad{\rm if}\quad \Re(\s)>0.
\ec\\
&\bc
&\hskip -0.3cm>\prod_{i=1}^g|(1-\om_i q^{2\s})(1-\ov\om_i q^{2\s})|\cdot |(q^{-2\s}-1)|\qquad{\rm if}\quad \Re(\s)<0\\
&\hskip -0.3cm <\prod_{i=1}^g|(1-\om_i q^{2\s})(1-\ov\om_i q^{2\s})|\cdot|(q^{-2\s}-1)|\qquad{\rm if}\quad \Re(\s)>0.
\ec
\ea
$$
which is nothing but the right hand side of \eqref{eq12}, provided that $g\geq 1$. This implies that unless $\Re(\s)=0$, i.e. unless $\s$ is a pure imaginary complex number, $\wh\zeta_{X,\F_q}^{\,\SL_2}(s)$ cannot be zero. This then proves the following
\begin{thm} (Yoshida) Let $X$ be an integral  regular projective  curve of genus $g$, then the $\SL_2$
zeta function
 $\wh\zeta_{X,\F_q}^{\,\SL_2}(s)$ and hence  the rank two zeta function $\wh\zeta_{X/\F_q,2}(s)$  satisfy the Riemann hypothesis.
\end{thm}
This result is due to H. Yoshida \cite{Y}, during our intensive lectures on non-abelian zeta functions for number fields in Kyoto. Yoshida, motivated by our works on rank two zeta functions of number fields \cite{Ki,LS,W2}, actually proves the Riemann hypothesis for a slight more general zeta function
\be\label{eq14}
\wh\zeta_{X,\F_q}(s;C_1):=C_1(s)\frac{\wh\zeta^{~}_{X/\F_q}(2s)}{1-q^{1-s}}-C_2(s)\frac{q^{-s}\,\wh\zeta^{~}_{X/\F_q}(2s-1)}{1-q^{2s-1}},
\ee
where $C_1(s)$ takes the form
\be
C_1(s)=q^{as}(1+q^{-s})q^{-hs}\prod_{j=1}^h(1-\ga_i q^{s-1/2})(1-\de_i q^{s-1/2})
\ee
with constants $\ga_i$ and $\de_i$ satisfying the conditions that $\ga_i+\de_i\in \R,\, |\ga_i+\de_i|\leq q+1$ for a non-negative real number $a$ and a natural number $h$, and $C_2(s)$ determined by
 \be
C_2(s)=C_1(1-s)\ee so that  \be
\wh\zeta_{X,\F_q}^{\,\SL_2}(1-s)=\wh\zeta_{X,\F_q}^{\,\SL_2}(s).\ee
Since our  $\SL_2$-zeta function $\wh\zeta_{X,\F_q}^{\,\SL_2}(s)$ is certainly a special form of Yoshida's type, the theorem above can be viewed as a special case.

\section{Applications of the RH for higher rank zetas}

\subsection{First estimates}

Before going further, we here deduce some natural upper and lower bounds for the non-abelian geo-arithmetic invariants $\a_{X,\F_q;n}(mn)\ (m=0,\ldots, g-1)$ and $\b_{X,\F_q;n}(0)$ of a curve $X/\F_q$ by assuming the Riemann hypothesis for rank $n$ zeta function of $X/\F_q$. Indeed as we will see later,   these bounds in lower ranks, in turn,  play key roles in proving the Riemann Hypothesis for higher rank  zeta functions.

Set now
\be
\a_{X,\F_q;n}'(mn):=\frac{\a_{X,\F_q;n}(mn)}{\a_{X,\F_q;n}(0)}\qqan \b_{X,\F_q;n}'(0)=\frac{\b_{X,\F_q;n}(0)}{\a_{X,\F_q;n}(0)}.\footnote{This definition make sense, since $\a_{X,\F_q;n}(0)\geq \frac{q^{h^0(X,\O_X^{\oplus n})}-1}{\Aut\,\O_X^{\oplus n}}>0$. In particular, $\a_{X/\F_q}(0)=\frac{q^{h^0(X,\O_X)}-1}{\#\Aut\,\O_X}=1$.}
\ee
Assume the Riemann hypothesis for the rank $n$ non-abelian zeta function of $X/\F_q$. By \eqref{eq5}, we have
$$
\ba
&\frac{1}{\a_{X,\F_q;n}(0)}P_{X,\F_q;n}(s)={T^{g-1}(1-T)(1-QT)}\cdot \frac{\wh Z_{X,\F_q;n}(T)}{\a_{X,\F_q;n}(0)}\\
%=&T^{g-1}\sum_{m=0}^{g-2}\a_{X,\F_q;n}(2m)\left(T^{m-(g-1)}+(QT)^{g-1-m}\right)+\a_{X,\F_q;n}(2(g-1))+\frac{(Q-1)T\beta_{X,\F_q;n}(0)}{(1-T)(1-QT)}\\
=&\Biggl(\left(\sum_{m=0}^{g-2}\a_{X,\F_q;n}'(mn)\left(T^{m}+(QT)^{2(g-1)-m}\right)+\a_{X,\F_q;n}'(n(g-1))T^{g-1}\right)(1-(Q+1)T+QT^2)\\
&\hskip 9.0cm+(Q-1)T^g\beta_{X,\F_q;n}'(0)\Biggr)\\
=&\prod_{i=1}^g(1-\om^{~}_{X/\F_q,n;i}T)(1-\ov \om^{~}_{X/\F_q,n;i}T).
\ea
$$
where $\om^{~}_{X/\F_q,n;i}$ denotes the reciprocal roots of $P_{X,\F_q;n}(s)$. Indeed,  $\frac{1}{\a_{X,\F_q;n}(0)}P_{X,\F_q;n}(s)\in \R[T]$ is a degree $2g$ polynomial in $T$ of real coefficients with leading coefficient $Q^g$ and constant term 1. By the functional equation of $\wh Z_{X,\F_q;n}(T)$, we may regroup all $2g$ reciprocal roots of $P_{X,\F_q;n}(s)$ into $g$ pairs, within each pair of which the products of two elements are always equal to $Q$. Consequently, the Riemann hypothesis for rank $n$ zeta function of $X/\F_q$ is equivalent to the condition that each such a  pair is of the form $\left\{\om^{~}_{X/\F_q,n;i},\,\ov\om^{~}_{X/\F_q,n;i}\right\}$. This is certainly equivalent to the condition that \be\label{eq41}
\left|\om^{~}_{X/\F_q,n;i}\right|=Q^\frac{1}{2}\qquad (i=1,\ldots g)
\ee
since $\om^{~}_{X/\F_q,n;i}\cdot \ov\om^{~}_{X/\F_q,n;i}=Q$. 

Set now $a^{~}_{X,\F_q,n;i}=\om^{~}_{X/\F_q,n;i}+ \ov\om^{~}_{X/\F_q,n;i}$. From Vieta's theorem between the reciprocal roots and coefficients of polynomials, by comparing the coefficients of $T^i$ for $i=1,\ldots, g$ in both sides of the above identity, under the Riemann hypothesis \eqref{eq41}, we conclude that the follows hold:
$$\ba
\bc\Big|\a_{X,\F_q;n}'(n)-(Q+1)\a_{X,\F_q;n}'(0)\Big|=\left|\sum_{i=1}^ga_{n,i}\right|\leq 2g\sqrt Q,\\
%~\\
\Big|\a_{X,\F_q;n}'(mn)-(Q+1)\a_{X,\F_q;n}'((m-1)n)+Q\a_{X,\F_q;n}'((m-2)n)\Big|\hskip 1.70cm (2\leq m\leq g-1)\\
\hskip 4.30cm\leq \sum _{1\leq i_{1}<i_{2}<\cdots <i_{m}\leq 2g}Q^{m/2}=\binom{2g}{m}Q^{m/2}\\
%~\\
\Big|(Q-1)\beta_{X,\F_q;n}'(0)-(Q+1)\a_{X,\F_q;n}'((g-1)n)+Q\a_{X,\F_q;n}'((g-2)n)\Big|\\
\hskip 4.30cm\leq \sum _{1\leq i_{1}<i_{2}<\cdots <i_{g}\leq 2g}Q^{g/2}=\binom{2g}{g}Q^{g/2}
\ec
\ea$$
Now expand each absolute value inequality, say  $|\kappa_i|<c_i$ as $-c_i<\kappa_i<c_i$, first, and then add three consecutive relations together with their lower and upper bounds, we arrive at
$$\ba
\bc
-2g\sqrt Q\leq\a_{X,\F_q;n}'(n)-(Q+1)\a_{X,\F_q;n}'(0)\leq 2g\sqrt Q,\\
%~\\
-\sum_{k=3}^{m+1}Q^{k-3}\sum_{i=1}^k\binom{2g}{i-1}Q^{(i-1)/2}-2gQ^{m-1}\sqrt Q\\
\hskip 1.50cm\leq \a_{X,\F_q;n}'(mn)-\Big(Q^m+\ldots+Q^2+Q+1\Big)\a_{X,\F_q;n}'(0)\hskip 1.50cm(2\leq m\leq g-1)\\
\hskip 3.0cm \leq\sum_{k=3}^{m+1}Q^{k-3}\sum_{i=1}^k\binom{2g}{i-1}Q^{(i-1)/2}+2gQ^{m-1}\sqrt Q,\\
%~\\
-\sum_{k=3}^{g+1}Q^{k-3}\sum_{i=1}^k\binom{2g}{i-1}Q^{(i-1)/2}-2gQ^{g-1}\sqrt Q\\
\hskip 1.50cm\leq (Q-1)\beta_{X,\F_q;n}'(0)-\Big(Q^{g}+\ldots+Q^2+Q+1\Big)\a_{X,\F_q;n}'(0)\\
\hskip 3.0cm\leq\sum_{k=3}^{g+1}Q^{k-3}\sum_{i=1}^k\binom{2g}{i-1}Q^{(i-1)/2}+2gQ^{g-1}\sqrt Q.\ec\ea$$

But, by definition, $\a_{X,\F_q;n}'(0)=1$. Therefore, we have proved the following:

\begin{prop}[Rough Bounds of non-abelian invariants]\label{thm3.3}
Assume the Riemann Hypothesis for the rank $n$ zeta functions of a projective regular integral curve $X$ over $\F_q$ of genus $g$, we have for the invariants $\beta_{X,\F_q;n}'(0)$ and $\a_{X,\F_q;n}'(2m)\ (m=1,\ldots,g-1)$
$$\ba
\bc
(Q+1)-2g\sqrt Q\leq\a_{X,\F_q;n}'(n)\leq (Q+1)+2g\sqrt Q\\
\\
\Big(Q^{m}+\ldots+Q^2+Q+1\Big)-\sum_{k=3}^{m+1}Q^{k-3}\sum_{i=1}^k\binom{2g}{i-1}Q^{(i-1)/2}-2gQ^{m-1}\sqrt Q\\
\hskip 2.0cm\leq \a_{X,\F_q;n}'(mn)\hskip 6.80cm(2\leq m\leq g)\\
\hskip 2.0cm\leq\Big(Q^{m}+\ldots+Q^2+Q+1\Big)+\sum_{k=3}^{m+1}Q^{k-3}\sum_{i=1}^k\binom{2g}{i-1}Q^{(i-1)/2}+2gQ^{m-1}\sqrt Q\\
\\
\Big(Q^{g}+\ldots+Q^2+Q+1\Big)-\sum_{k=3}^{g+1}Q^{k-3}\sum_{i=1}^k\binom{2g}{i-1}Q^{(i-1)/2}-2gQ^{g-1}\sqrt Q\\
\hskip 2.0cm\leq (Q-1)\beta_{X,\F_q;n}'(0)\\
\hskip 2.0cm\leq\Big(Q^{g}+\ldots+Q^2+Q+1\Big)+\sum_{k=3}^{g+1}Q^{k-3}\sum_{i=1}^k\binom{2g}{i-1}Q^{(i-1)/2}+2gQ^{g-1}\sqrt Q
\ec
\ea$$ 
\end{prop}
Even somehow rough upper and lower bounds are obtained here for the $\a$ and $\b$-invariants in rank $n$, nevertheless, as we will see in \S\ref{sec3.5} below, these bounds are far from being sharp.

\subsection{Counting Brill-Noether stratifications}
In this subsection, we explain why the rank $n$ non-abelian $\a$ and $b$ invariants are very important, in order to understand the intrinsic structures of $X/\F_q$.

\begin{exam}[Brill-Noether Loci] Even when $n=1$, the above result exposes some intrinsic geo-arithmetic properties of the curve $X/\F_q$. Indeed, we may introduce the Brill-Noether loci within  the degree Picard group $\Pic^d(X)$ of $X$ by setting
\be
W_{X/\F_q}^{\geq i}(d):=\big\{L\in \Pic^d(X):h^0(X,L)\geq i\big\}\qan W_{X/\F_q}^{=i}(d):=\big\{L\in \Pic^d(X):h^0(X,L)= i\big\}.
\ee
The $W_{X/\F_q}^{\geq i}(d)$'s induce a natural stratification structure on $\Pic^d(X)$ since
\be
W_{X/\F_q}^{\geq i}(d)=\bigsqcup_{j\geq i}W_{X/\F_q}^{=j}(d)\qqan W_{X/\F_q}^{\geq 0}(d)=\Pic^d(X).
\ee
It is natural to ask what are the topological or better motivic properties of these refined structures. 
Set accordingly
\be
w_{X/\F_q}^{\geq i}(d)=\#W_{X/\F_q}^{\geq i}(d)\qqan w_{X/\F_q}^{=i}(d)=\#W_{X/\F_q}^{=i}(d).
\ee
and 
\be
\a_{X/\F_q}^{\geq i}(d)=\sum_{L\in W_{X/\F_q}^{\geq i}(d)}\frac{q^{h^0(X,L)-1}}{q-1}\qan \a_{X/\F_q}^{=i}(d)=\sum_{L\in W_{X/\F_q}^{= i}(d)}\frac{q^{h^0(X,L)-1}}{q-1}=\frac{q^i-1}{q-1}w_{X/\F_q}^{=i}(d).
\ee
Then  
\be
w_{X/\F_q}^{\geq 0}(d)=\sum_{i\geq 0}w_{X/\F_q}^{=i}(d)=\#\Pic^0(X)\qan\a_{X/\F_q}(d)=\a_{X/\F_q}^{\geq 0}(d)=\sum_{i\geq 0}\frac{q^i-1}{q-1}w_{X/\F_q}^{=i}(d),
\ee
which by Theorem\,\ref{thm3.3}, is controlled by
\be\label{eq46}
\bc
-2g\sqrt q\leq\a_{X/\F_q}(1)-(q+1)\leq 2g\sqrt q\\
%~\\
-\sum_{k=3}^{m+1}q^{k-3}\sum_{i=1}^k\binom{2g}{i-1}q^{(i-1)/2}-2gq^{m-1}\sqrt q\\
\hskip 2.0cm\leq \a_{X/\F_q}'(m)-\Big(q^{m}+\ldots+q^2+q+1\Big)\hskip 1.80cm(2\leq m\leq g)\\
\hskip 4.0cm\leq\sum_{k=3}^{m+1}q^{k-3}\sum_{i=1}^k\binom{2g}{i-1}q^{(i-1)/2}+2gq^{m-1}\sqrt q\\
%~\\
-\sum_{k=3}^{g+1}q^{k-3}\sum_{i=1}^k\binom{2g}{i-1}q^{(i-1)/2}-2gq^{g-1}\sqrt q\\
\hskip 2.0cm\leq (q-1)\beta_{X/\F_q}(0)-\Big(q^{g}+\ldots+q^2+q+1\Big)\\
\hskip 4.0cm\leq\sum_{k=3}^{g+1}q^{k-3}\sum_{i=1}^k\binom{2g}{i-1}q^{(i-1)/2}+2gq^{g-1}\sqrt q
\ec
\ee
since $\a_{X/\F_q}(0)=1$.

Recall that, by definition, 
\be
\beta_{X/\F_q}(0)=\frac{1}{q-1}\#\Pic^0(X)=\wh v_1=\frac{1}{q-1}\prod_{i=1}^g(1-\om_{X/F_q,i})(1-\ov\om_{X/F_q,i}).\ee
Hence it admits a natural bound
\be
\frac{1}{q-1}\prod_{i=1}^g(1-\sqrt q)(1-\sqrt q)\leq\beta_{X/\F_q}(0)\leq \frac{1}{q-1}\prod_{i=1}^g(1+\sqrt q)(1+\sqrt q)
\ee
So, at least for the $\b$-invariant, \eqref{eq46} may not be the sharpest bounds.

Furthermore, we may introduce natural associated invariants
\be
\a_{X/\F_q}^{\geq i}(d):=\sum_{V\in W_{X/\F_q}^{\geq i}(d)}\frac{q^{h^0(X,L)}-1}{q-1}
=\sum_{j\geq i}w_{X/\F_q}^{=j}(d)\frac{q^j-1}{q-1}
\ee
and their associated generating function, for fixed $i$ and $d$,
\be
A_{X/\F_q}(u,v):=\sum_{i,d\geq 0}a_{X/\F_q}^{=i}(d)\cdot u^iv^d.
\ee
Recall that, by the vanishing theorem, we have
\be
\a_{X/\F_q}^{=i}(d)=\delta_{i,d-(g-1)}\frac{q^{d-(g-1)}-1}{q-1}\#\Pic^0(X) \qquad d\geq g
\ee
So,
$A_{X/\F_q}(u,v)$ is indeed a rational function of $u$ and $v$. A modification can be given to result a rational function which  satisfies the  standard functional equation. The detailed discussion will be given elsewhere, together with the associated Riemann hypothesis property.

Obviously, these structures admit a natural generalization to the setting of rank $n$ semi-stable bundles. For instances,
we expect that the non-abelian motivic structures of $X/\F_q$ may be understood via the rank $n$ generating function
\be
B(X,\F_q;u,v;z):=\sum_{n\geq 0}A(X,\F_q;n;u,v)z^n=\sum_{n,i,m\geq 0}\a_{X,\F_q;n}^{=i}(mn)u^iv^{mn}z^n.
\ee
Similarly, by the vanishing theorem, we have
\be
\a_{X/\F_q}^{=i}(mn)=\delta_{i,mn-n(g-1)}\frac{q^{d-n(g-1)}-1}{q-1}\#M_{X,\F_q;n}(0) \qquad d\geq ng
\ee
where $M_{X,\F_q;n}(0)$ denotes the space of isomorphism classes, or better by including the
geo-arithmetic structures, the moduli stack, of semi-stable vector bundles
of rank $n$ and degree $mn$ on $X$ rationally over $\F_q$.
So,
$B_{X/\F_q}(u,v;z)$ is indeed a rational function of $u$ and $v$, for which a modification can be made to result a rational function that satisfies the standard functional equation. For details, please refer to our work in preparation on motivic structures of curves over general base fields.
\end{exam}
We end this discussion with the following application of the Riemann hypothesis.
\begin{lem}\label{cor3.2}
For a regular projective  curve $X$ of genus $g\geq 2$ over $\F_q$, we have
\be\label{eq58}
\left(q^{-\binom{n}{2}(g-1)}\cdot\frac{\b_{X,\F_q;2}(0)}{\wh\nu_1^2}-\frac{1}{q+1}\right)>0.
\qquad
{\rm In\ particular,}\qquad
\frac{\wh v_{2}}{\wh v_1^{\,2}}>\frac{q}{q^2-1}
\ee
Slight differently,
\be\label{eq60}
\frac{\wh v_{2}}{\wh v_1^{\,2}}>\frac{3}{2}\frac{1}{q-1}.
\ee
\end{lem}
\bp
Note that
$$
\left(q^{-\binom{n}{2}(g-1)}\cdot\frac{\b_{X,\F_q;2}(0)}{\wh\nu_1^2}-\frac{1}{q+1}\right)
=\frac{\wh\nu_2+\frac{1}{1-q^2}\wh\nu_1^2}{\wh\nu_1^2}-\frac{1}{q+1}
%=\frac{\wh\nu_2}{\wh\nu_1^2}+\frac{1}{1-q^2}-\frac{1}{q+1}\\
%=&\frac{\wh\nu_2}{\wh\nu_1^2}+\frac{1-1+q}{1-q^2}
=\frac{\wh\nu_2}{\wh\nu_1^2}-\frac{q}{q^2-1}
$$
Now
$$\ba
\frac{\wh\nu_2}{\wh\nu_1^2}=&\left(\prod_{i=1}^g\frac{(1-\om_{X/\F_q;i}q^{-2})(1-\ov\om_{X/\F_q;i}q^{-2})}{(1-\om_{X/\F_q;i}q^{-1})^2(1-\ov\om_{X/\F_q;i}q^{-1})^2}\right)\frac{q^{-2(g-1)}(q-1)^2}
{q^{-2(g-1)}(1-q^-2)(1-qq^{-2})}\\
%=&\left(\prod_{i=1}^g\frac{(1-\om_{X/\F_q;i}q^{-2})(1-\ov\om_{X/\F_q;i}q^{-2})}{(1-\om_{X/\F_q;i}q^{-1})^2(1-\ov\om_{X/\F_q;i}q^{-1})^2}\right)\frac{q^3}{q+1}\\
%=&\left(\prod_{i=1}^g\frac{(1-\om_{X/\F_q;i}q^{-2})(1-\ov\om_{X/\F_q;i}q^{-2})}{(1-\om_{X/\F_q;i}q^{-1})^2(1-\ov\om_{X/\F_q;i}q^{-1})^2}\right)\frac{q^3}{q+1}\\
%=&\left(\prod_{i=1}^g\frac{(1-\ov\om_{X/\F_q;i}^{-1}q^{-1})(1-\om_{X/\F_q;i}^{-1}q^{-1})}{(1-\om_{X/\F_q;i}q^{-1})(1-\ov\om_{X/\F_q;i}q^{-1})}\frac{1}{(1-\om_{X/\F_q;i}q^{-1})(1-\ov\om_{X/\F_q;i}q^{-1})}\right)\frac{q^3}{q+1}\\
=&\frac{1}{q^g}\left(\prod_{i=1}^g\frac{(\ov\om_{X/\F_q;i}-q^{-1})(\om_{X/\F_q;i}-q^{-1})}{(1-\om_{X/\F_q;i}q^{-1})(1-\ov\om_{X/\F_q;i}q^{-1})}\frac{1}{(1-\om_{X/\F_q;i}q^{-1})(1-\ov\om_{X/\F_q;i}q^{-1})}\right)\frac{q^3}{q+1}\\
&\hskip 7.0cm({\rm since}\ \om_{X/\F_q;i}\ov \om_{X/\F_q;i}=q)\\
%=&\prod_{i=1}^g\left(\frac{(\ov\om_{X/\F_q;i}-q^{-1})(\om_{X/\F_q;i}-q^{-1})}{(1-\om_{X/\F_q;i}q^{-1})(1-\ov\om_{X/\F_q;i}q^{-1})}\frac{q^{-(g-1)}(q-1)}{(1-\om_{X/\F_q;i}q^{-1})(1-\ov\om_{X/\F_q;i}q^{-1})}\right)\frac{q^2}{q^2-1}\\
=&\left(\prod_{i=1}^g\frac{(\ov\om_{X/\F_q;i}-q^{-1})(\om_{X/\F_q;i}-q^{-1})}{(1-\om_{X/\F_q;i}q^{-1})(1-\ov\om_{X/\F_q;i}q^{-1})}\right)\frac{q^{2}}
{q^2-1}\cdot \wh \nu_1
\geq \left(\prod_{i=1}^g1\right)\frac{q^{2}}
{q^2-1}\cdot \wh\nu_1\\
&\hskip 2.0cm{\rm by\ Yoshida's\ lemma,\
since\ }|q^{-1}|<1\ {\rm (and\ }\ \wh\nu_1>0\ {\rm and}\ \frac{\wh\nu_2}{\wh\nu_1^2}>0)\\
=&\frac{q^{2}}
{q^2-1}\cdot\wh\nu_1=\frac{q^{2}}
{(q^2-1)(q-1)}\cdot\#\Pic^0(X)>\frac{q}{q^2-1}\\
\ea$$
where in the last step, we have used the fact that there exists at least one $\F_q$ rational point in 
$\Pic^0(X)$.
This verifies \eqref{eq58}. 

In addition, since $\#\Pic^0(X)\geq 2$, we have $2q^2\#\Pic^0(X)+3>3q^2$. This implies that
$\frac{q^2}{q^2-1}\#\Pic^0(X)>\frac{3}{2}$. Therefore,
\be
\frac{\wh v_{2}}{\wh v_1^{\,2}}>\frac{q^2}{(q^2-1)(q-1)}\#\Pic^0(X)>\frac{3}{2}\frac{1}{q-1},
\ee
as wanted.
\ep
This result will be used in the next section to prove the rank three Riemann hypothesis. In addition, an alternative, but much neat proof will be provided in the following subsection.

\subsection{Rank $n$ beta invariant $\b_{X,\F_q;n}(0)$}

In this subsection, we examine in details the structures of the rank $n$ beta invariant $\b_{X,\F_q;n}(0)$ for an integral regular projective curve $X$ over a finite field $\F_q$.

Let $X$ be an integral regular projective curve of genus $g$ over a finite field $\F_q$ of $q$ elements.
For a fixed positive integer $n\in\Z_{>0}$,  the rank $n$ non-abelian beta invariant $\b_{X,\F_q;n}(d)$ is defined as
\be
\b_{X,\F_q;n}(d)=\sum_{\cE}\frac{1}{\#\Aut(\cE)}
\ee
where $\cE$ runs through all semi-stable $\F_q$-rational vector bundles of rank $n$ and degree $d$ on $X$. They are introduced by Harder-Narasimhan \cite{HN}, in their study of Betti numbers of the moduli spaces of semi-stable vector bundles on $X$.

As a direct consequence of Siegel's formula on the total mass of rank $n$ vector bundles of degree $d$ on $X$ with fixed determinants associated to the so-called Tamagawa measure, if we introduce $\b_{X,\F_q;n}^\tot$ instead by taking the sum on all rank $n$ vector bundles of degree $d$ on $X/\F_q$, we have,  see e.g. \cite{DR},
\be
\b_{X,\F_q;n}^\tot=q^{\binom{n}{2}(g-1)}\cdot \wh\nu_{X/\F_q}^{~}(n):=q^{\binom{n}{2}(g-1)}\prod_{k=1}^n\wh\zeta_{X/\F_q}(k).
\ee 
Here 
\be
\wh\nu_{X/\F_q}^{~}(n):=\prod_{k=1}^n\wh\zeta_{X/\F_q}(k),\qquad
\wh\zeta_{X/\F_q}(1):=\Res_{t=1}\wh Z_{X/\F_q}(t)\ee 
and $\wh\zeta_{X/\F_q}(k)\ (k\geq 2)$
denotes the special value at $k$
 of the complete Artin zeta/Zeta function  of $X/\F_q$
\be
\wh\zeta_{X/\F_q}(s):=\sum_{D\geq 0}\bl q^{\chi(D)}\br^{-s}:=\wh Z_{X/\F_q}(t)\qquad (t:=q^{-s}),
\ee
where, as usual, $D$ run through all effective divisors of $X$. 
Furthermore, to understand the $\b$-invariants $\b_{X,\F_q;n}(d)$, Harder-Narasimhan introduces
parabolic reduction \cite{HN}. The up-shot of their method is  the following inductive formula for $\b_{X,\F_q;n}(d)$:

\be
\beta_{X,\F_q;n}^\tot(d)=\sum_{\substack{n_1,\ldots,n_k\geq 1\\ n_1+\ldots+n_k=n}}\nu_{X/\F_q}^{~}(1)^{k-1}\sum_{\substack{(\de_i)\in\prod_{i=1}^k\Z/n_i\\ d_i\equiv\de_i(n_i)\\i=1,\ldots,k}}
\sum_{\substack{d_1,\ldots,d_k\\ d_1+\ldots+d_k=d\\ \frac{d_1}{n_1}>\ldots>\frac{d_k}{n_k}
}}q^{\sum_{i<j}\bl n_in_j(g-1)-(n_id_j-n_jd_i)\br}\prod_{i=1}^k\b_{X,\F_q;n_i}(\de_i)\ee 
where $\nu_{X/\F_q}(1):=\Res_{t=1}\zeta_{X/\F_q}(t)$ denotes the residue at $t=1$ of the Artin Zeta function 
\be
Z_{X/\F_q}(t):=\zeta_{X/\F_q}(s):=\sum_{D\geq 0}\bl q^{\deg(D)}\br^{-s}.\ee
In particular,
$$
\ba
\beta_{X,\F_q;n}(d)\=&\beta_{X,\F_q;n}^\tot(d)\\
-&\sum_{\substack{k\geq 2\\n_1,\ldots,n_k\geq 1\\ n_1+\ldots+n_k=n}}\nu_{X/\F_q}^{~}(1)^{k-1}\sum_{\substack{(\de_i)\in\prod_{i=1}^k\Z/n_i\\ d_i\equiv\de_i(n_i)\\i=1,\ldots,k}}
\sum_{\substack{d_1,\ldots,d_k\\ d_1+\ldots+d_k=d\\ \frac{d_1}{n_1}>\ldots>\frac{d_k}{n_k}
}}q^{\sum_{i<j}\bl n_in_j(g-1)-(n_id_j-n_jd_i)\br}\prod_{i=1}^k\b_{X,\F_q;n_i}(\de_i)\ea$$
These formulas are finally simplified by Zagier \cite{Z} to obtain a finite closed form
\begin{thm} Let $X$ be an integral regular projective curve of genus $g$ ovcer a finite field $\F_q$. Then its rank $n$ beta invariant $\b_{X,\F_q;n}(d)$ is given by
\be
\beta_{X,\F_q;n}(d)=q^{\binom{n}{2}(g-1)}\cdot\sum_{\substack{n_1,\ldots,n_k\geq 1\\ n_1+\ldots+n_k=n}}\frac{q^{(n_i+n_{i+1})\{(n_1+\ldots+n_i)\frac{d}{n}\}}}{\prod_{i=1}^{k-1}(1-q^{n_i+n_{i+1}})}\prod_{i=1}^k\wh \nu_{X/\F_q}(n_i).
\ee
In particular, when $n|d$, say $d=mn$ for a certain $m\in\Z$, we have
\be
\beta_{X,\F_q;n}(mn)=q^{\binom{n}{2}(g-1)}\cdot\sum_{\substack{n_1,\ldots,n_k\geq 1\\ n_1+\ldots+n_k=n}}\frac{\prod_{i=1}^k\wh \nu_{X/\F_q}(n_i)}{\prod_{i=1}^{k-1}(1-q^{n_i+n_{i+1}})}.
\ee
\end{thm}
Note  that, induced from the map $\cE\mapsto \cE\otimes \O_C(A)$ for a fixed degree one $\F_q$-rational divisor $A$ on $X$, whose existence is a well-known consequence of the theory
of Artin zeta functions of curves over finite fields (see e.g. \cite{Me}), we arrive at an elementary relation
 \be
\beta_{X,\F_q;n}(mn)=\beta_{X,\F_q;n}(0)\qquad(\forall m\in\Z),
\ee
Besides this, in general, it is quite complicated to estimate
the $\b$-invariants for a curve $X$, and hence these $\b$-invariants  can hardly be controlled. 

However, after the introduction of the rank $n$ non-abelian zeta function $\wh\zeta_{X,\F_q;n}(s)$ of $X/\F_q$ in \cite{W}, the situation has changed dramatically.
Recall that, by definition,
\be
\wh\zeta_{X,\F_q;n}(s):=\sum_{\cE}\frac{q^{h^0(X,\cE)}-1}{\Aut(\cE)}\big(q^{-s}\big)^{\chi(X,\cE)}.
\ee
where $\cE$ runs through all semi-stable $\F_q$-rational vector bundles of rank $n$ and degree $mn \ (m\geq 0)$ on $X$.
Moreover,  as said earlier, by the Riemann-Roch theorem, cohomology duality and the vanishing theorem for semi-stable vector bundles, tautologically, using the standard zeta technique, we have,
 for $Q:=q^n$ and $T:=t^n$ with $t=q^{-s}$,
$$\ba
&\wh\zeta_{X,\F_q;n}(s)=\wh Z_{X,\F_q;n}(s)\\
=&\sum_{m=0}^{g-2}\a_{X,\F_q;n}(mn)\bl T^{m-(g-1)}+(QT)^{(g-1)-m}\br\br+\a_{X,\F_q;n}(n(g-1))+\frac{(Q-1)\beta_{X,\F_q;n}(0)\cdot T}{(1-T)(1-QT)}
\ea
$$
In particular,  $\wh\zeta_{X,\F_q;n}(s)$ is a rational function of $T$, satisfies the standard functional equation $s\leftrightarrow 1-s$
and admits only two singularities as a function of $T$, namely, simple poles at $T=1$ and $T=Q^{-1}$ with residue
\be
\Res_{T=Q^{-1}}\wh\zeta_{X,\F_q;n}(s)=\b_{X,\F_q;n}(0).
\ee
In addition, when $n=1$, $\wh\zeta_{X,\F_q;n}(s)$ coincides with the completed Artin zeta function $\wh\zeta_{X/\F_q}(s)$ of $X/\F_q$. 
 In this way, the $\b$-invariant $\b_{X,\F_q;n}(0)$ in rank $n$ of $X$ can be understood  in terms of the theory of non-abelian zeta functions of curves. In particular, we expect that the Riemann hypothesis holds for these non-abelian zeta functions as well: All zeros of the rank $n$-non-abelian zeta function $\wh\zeta_{X,\F_q;n}(s)$ of $X/\F_q$ lie on the central line $\Re(s)=\frac{1}{2}$.

 In the sequel, till the end of this subsection,  we assume that this Riemann hypothesis  holds for all $\wh\zeta_{X,\F_q;k}(s)$ with  $k\leq n$. Accordingly, we have
\be\label{eq13}
\wh\zeta_{X,\F_q;n}(s)=\wh Z_{X,\F_q;n}(T)=\a_{X,\F_q;n}(0)\cdot \frac{\prod_{i=1}^g(1-\om_{X,\F_q;n}(i)T)(1-\ov \om_{X,\F_q;n}(i)T)}{(1-T)(1-QT)T^{g-1}}
\ee
where $\om_{X,\F_q;n}(i)$ denotes the $i$-th reciprocal root of $\wh Z_{X,\F_q;n}(T)$. Recall that by the functional equation, we have
\be
\om_{X,\F_q;n}(i)\cdot \ov \om_{X,\F_q;n}(i)=Q\qquad (i=1,\ldots,g)
\ee
Furthermore, from our assumption that the Riemann Hypothesis in rank $n$ holds, we have 
\be
\big|\om_{X,\F_q;n}(i)\big|=\sqrt Q\qquad (i=1,\ldots,g)
\ee
Therefore, by \eqref{eq13}, we have
$$
\ba
\Res^{~}_{T=1}\wh\zeta_{X,\F_q;n}(s)=&\a_{X,\F_q;n}(0)\cdot \frac{\prod_{i=1}^g(1-\om_{X,\F_q;n}(i))(1-\ov \om_{X,\F_q;n}(i))}{(Q-1)}\\
=&\a_{X,\F_q;n}(0)\cdot \frac{\prod_{i=1}^g\bl 1-a_{X,\F_q;n}(i)+Q\br}{(Q-1)}
\ea
$$
where $
a_{X,\F_q;n}(i):=\om_{X,\F_q;n}(i)+\ov \om_{X,\F_q;n}(i)$ satisfying 
\be |a_{X,\F_q;n}(i)|\leq 2\sqrt Q\qquad (i=1,\ldots,g).\ee

In this way, we have proved the following
\begin{cor}\label{cor3.3}  Let $X$ be an integral  regular projective curve of genus $g\geq 2$ over a finite field $\F_q$. Assume that the rank $n$ Zeta function $\wh Z_{X,\F_q;n}(T)$ of $X/\F_q$ satisfies the Riemann hypothesis. Then
\be
\frac{\ \bl \sqrt Q-1\br^{2g-1}\ }{\sqrt Q+1}\leq \b_{X,\F_q;n}'(0)=\frac{\ \prod_{i=1}^g\bl 1-a_{X,\F_q;n}(i)+Q\br\ }{(Q-1)}\leq
\frac{\ \bl \sqrt Q+1\br^{2g-1}\ }{\sqrt Q-1}
\ee
\end{cor}

\subsection{An application  of the counting miracle}

Recall that 
\be
\b_{X,\F_q;n}'(0)=\frac{\b_{X,\F_q;n}(0)}{\a_{X,\F_q;n}(0)}.
\ee
In this way, it suffices to understand the $\a$-invariants $\a_{X,\F_q;n}(0)$. For this, let  $C_{n,g;q}:=~q^{\binom{n}{2}(g-1)}$. By applying
the first relation to $m=0$ of Theorem\,\ref{thm1.7}, note that 
\be
\a_{X/\F_q}(0)=1
\ee and that the first summation is understood as 0 when $g=1$, we obtain
 $$\ba
&\frac{1}{C_{n,g;q}}\a_{X,\F_q;n}(0)
%=&\sum_{a=1}^nq^{(n-a)(0-(g-1))}\Biggl(
%-\sum_{\substack{k_1,\ldots,k_p>0\\ k_1+\ldots+k_p=n-a}}\frac{\wh v_{k_1}\ldots\wh v_{k_p}}{\prod_{j=1}^{p-1}(1-q^{k_j+k_{j+1}})} \sum_{\substack{l_1,\ldots,l_r>0\\ l_1+\ldots+l_r=a-1}}
% \frac{\wh v_{l_1}\ldots\wh v_{l_r}}{\prod_{j=1}^{r-1}(1-q^{l_j+l_{j+1}})}\\
 %&\hskip 5.0cm\times\sum_{0=k+\ell+\kappa}\sum_{k=0}^{g-2}\a_{X/\F_q}(k)\sum_{\ell=1}^\infty\bl q^{-k_{p}}\br^\ell\sum_{\kappa=0}^\infty\bl q^{1+l_{1}}\br^\kappa\Biggr)\\
%%&\quad+\delta_{0,g-1}\Biggl(\a_{X/\F_q}\bl (g-1)\br\cdot\sum_{\substack{l_1,\ldots,l_r>0\\ l_1+\ldots+l_r=n-1}}
% \frac{\wh v_{l_1}\ldots\wh v_{l_r}}{\prod_{j=1}^{r-1}(1-q^{l_j+l_{j+1}})}\Biggr)\\
% =&\sum_{a=n}q^{(n-a)(0-(g-1))}\Biggl( \sum_{\substack{l_1,\ldots,l_r>0\\ l_1+\ldots+l_r=a-1}}
% \frac{\wh v_{l_1}\ldots\wh v_{l_r}}{\prod_{j=1}^{r-1}(1-q^{l_j+l_{j+1}})}\sum_{0=k+\kappa}\sum_{k=0}\a_{X/\F_q}(k)\sum_{\kappa=0}\bl q^{1+l_{1}}\br^\kappa\Biggr)\\
%&\quad+\delta_{0,g-1}\Biggl(\a_{X/\F_q}\bl (g-1)\br\cdot\sum_{\substack{l_1,\ldots,l_r>0\\ l_1+\ldots+l_r=n-1}}
 %\frac{\wh v_{l_1}\ldots\wh v_{l_r}}{\prod_{j=1}^{r-1}(1-q^{l_j+l_{j+1}})}\Biggr)\\
%  =&q^{(n-n)(0-(g-1))}\Biggl( \sum_{\substack{l_1,\ldots,l_r>0\\ l_1+\ldots+l_r=n-1}}
%% \frac{\wh v_{l_1}\ldots\wh v_{l_r}}{\prod_{j=1}^{r-1}(1-q^{l_j+l_{j+1}})}\a_{X/\F_q}(0)\Biggr)\\
% &\qquad+\delta_{g,1}\Biggl(\a_{X/\F_q}\bl (0)\br\cdot\sum_{\substack{l_1,\ldots,l_r>0\\ l_1+%\ldots+l_r=n-1}}
% \frac{\wh v_{l_1}\ldots\wh v_{l_r}}{\prod_{j=1}^{r-1}(1-q^{l_j+l_{j+1}})}\Biggr)\\
 =\sum_{\substack{l_1,\ldots,l_r>0\\ l_1+\ldots+l_r=n-1}}
 \frac{\wh v_{l_1}\ldots\wh v_{l_r}}{\prod_{j=1}^{r-1}(1-q^{l_j+l_{j+1}})}.\\
\ea$$

Indeed, as said earlier, this is a consequence of the special uniformity of zetas for curves over finite fields, claiming that, up to a non-zero constant multiple $C_{n,g;q}$, 
$$
 \ba
 \wh\zeta_{X,\F_q;n}(s)=&\wh\zeta_{X,\F_q}^{\SL_n}(s)=C_{n,g;q}\cdot \sum_{a=1}^{n}\sum_{\substack{k_1,\ldots,k_p>0\\ k_1+\ldots+k_p=n-a}}\frac{\wh v_{k_1}\ldots\wh v_{k_p}}{\prod_{j=1}^{p-1}(1-q^{k_j+k_{j+1})}} \frac{1}{(1-q^{ns-n+a+k_{p}})}\\
 &\hskip 2.0cm\times\wh\zeta_{X,\F_q}(ns-n+a)\sum_{\substack{l_1,\ldots,l_r>0\\ l_1+\ldots+l_r=a-1}}
 \frac{1}{(1-q^{-ns+n-a+1+l_{1}})}\frac{\wh v_{l_1}\ldots\wh v_{l_r}}{\prod_{j=1}^{r-1}(1-q^{l_j+l_{j+1})}}.
 \ea
 $$
Furthermore,
$$\ba
\b_{X,\F_q;n}(0)=&\Res_{T=1}\wh Z_{X,\F_q;n}\\
=&\Res_{T=1}\Biggl(C_{n,g;q}\cdot \sum_{a=1}^{n}\sum_{\substack{k_1,\ldots,k_p>0\\ k_1+\ldots+k_p=n-a}}\frac{\wh v_{k_1}\ldots\wh v_{k_p}}{\prod_{j=1}^{p-1}(1-q^{k_j+k_{j+1})}} \frac{1}{(1-q^{ns-n+a+k_{p}})}\\
 &\times\wh\zeta_{X,\F_q}(ns-n+a)\sum_{\substack{l_1,\ldots,l_r>0\\ l_1+\ldots+l_r=a-1}}
 \frac{1}{(1-q^{-ns+n-a+1+l_{1}})}\frac{\wh v_{l_1}\ldots\wh v_{l_r}}{\prod_{j=1}^{r-1}(1-q^{l_j+l_{j+1})}}\Biggr)\\
%=&\Res_{T=1}\Biggl(C_{n,g;q}\cdot \sum_{a=1}^{n}\sum_{\substack{k_1,\ldots,k_p>0\\ k_1+\ldots+k_p=n-a}}\frac{\wh v_{k_1}\ldots\wh v_{k_p}}{\prod_{j=1}^{p-1}(1-q^{k_j+k_{j+1})}} \frac{T}{(T-q^{-n+a+k_{p}})}\\
% &\times\frac{\prod_{i=1}^g(1-a_{X/\F_q}(i)Tq^{n-a}+Q(q^{n-a}T)^2)}{(1-Tq^{n-a})(1-QTq^{n-a})(q^{n-a}T)^{g-1}}\sum_{\substack{l_1,\ldots,l_r>0\\ l_1+\ldots+l_r=a-1}}
% \frac{1}{(1-q^{n-a+1+l_{1}}T)}\frac{\wh v_{l_1}\ldots\wh v_{l_r}}{\prod_{j=1}^{r-1}(1-q^{l_j+l_{j+1})}}\Biggr)\\
%=&\Res_{T=1}\Biggl(C_{n,g;q}\cdot \sum_{a=1}^{n}\sum_{\substack{k_1,\ldots,k_p>0\\ k_1+\ldots+k_p=n-a}}\frac{\wh v_{k_1}\ldots\wh v_{k_p}}{\prod_{j=1}^{p-1}(1-q^{k_j+k_{j+1})}} \frac{1}{(T-q^{-n+a+k_{p}})}\\
% &\times\frac{\prod_{i=1}^g(1-a_{X/\F_q}(i)q^{n-a}+Q(q^{n-a})^2)}{(1-Tq^{n-a})(1-QTq^{n-a})(q^{n-a}T)^{g-1}}\sum_{\substack{l_1,\ldots,l_r>0\\ l_1+\ldots+l_r=a-1}}
% \frac{1}{(1-q^{n-a+1+l_{1}}T)}\frac{\wh v_{l_1}\ldots\wh v_{l_r}}{\prod_{j=1}^{r-1}(1-q^{l_j+l_{j+1})}}\Biggr)\\
\ea$$
From here, proceeds as what has been done in the final section of \cite{WZ2}, we get 
\be
\b_{X,\F_q;n}(0)=C_{n,g;q}\sum_{\substack{k_1,\ldots,k_p>0\\ k_1+\ldots+k_p=n}}\frac{\wh v_{k_1}\ldots\wh v_{k_p}}{\prod_{j=1}^{p-1}(1-q^{k_j+k_{j+1})}} 
\ee
Therefore, 
\be
\a_{X,\F_q,n}(0)
=C_{n,g;q}\cdot\frac{\b_{X,\F_q;n-1}(0)}{C_{n-1,g}}=\frac{C_{n,g;q}}{C_{n-1,g}}\cdot \b_{X,\F_q;n-1}(0)
\ee
In particular, in order to make this relation coincide with that of Sugahara stated in the theorem of the appendix to \S1.4, we should have
$C_{n,g;q}=q^{\binom{n}{2}(g-1)}$, as claimed. This is compatible with the main theorem of \cite{MR} as well (even not with the formula (2.9) there). In this way, we have obtained yet another totally different proof of the following
\begin{thm}[Counting Miracle]\label{thm3.5} Let $X$ be an integral  regular projective curve of genus $g$ over a finite field $\F_q$.  Then we have, with $C_{n,g;q}:=q^{\binom{n}{2}(g-1)}$,
$$
\ba
\a_{X,\F_q;n}(0)=C_{n,g;q}\sum_{\substack{k_1,\ldots,k_p>0\\ k_1+\ldots+k_p=n-1}}\frac{\wh v_{k_1}\ldots\wh v_{k_p}}{\prod_{j=1}^{p-1}(1-q^{k_j+k_{j+1})}} \\
\b_{X,\F_q,n-1}(0)=C_{n-1,g}\sum_{\substack{k_1,\ldots,k_p>0\\ k_1+\ldots+k_p=n-1}}\frac{\wh v_{k_1}\ldots\wh v_{k_p}}{\prod_{j=1}^{p-1}(1-q^{k_j+k_{j+1})}} \\
\ea
$$
In particular, 
\be
\a_{X,\F_q,n+1}(0)
=q^{\left(\binom{n+1}{2}-\binom{n}{2}\right)(g-1)}\cdot \b_{X,\F_q;n}(0)=q^{n(g-1)}\cdot \b_{X,\F_q;n}(0)
\ee
\end{thm}

\subsection{Second estimates: upper and lower bounds for $\b$-invariants}\label{sec3.5}

Now, by Corollary\,\ref{cor3.3} and Theorem\,\ref{thm3.5}, we conclude that
\be
\frac{\b_{X,\F_q;n}(0)}{\a_{X,\F_q;n}(0)}=\frac{\b_{X,\F_q;n}(0)}{q^{(n-1)(g-1)}\b_{X,\F_q;n-1}(0)}=\frac{\ \prod_{i=1}^g\bl 1-a_{X,\F_q;n}(i)+q^n\br\ }{(q^n-1)}
\ee
Therefore, under the assumption that the Riemann hypothesis holds for all $\wh Z_{X,\F_q;k}(T_k)$ where $T_k=t^k, k=2,\ldots, n$, we have
$$\ba
\b_{X,\F_q;n}(0)=&{\b_{X,\F_q;1}(0)}\cdot\prod_{k=2}^{n}\frac{\b_{X,\F_q;k}(0)}{\b_{X,\F_q;k-1}(0)}
\\
=&{\b_{X,\F_q;1}(0)}\cdot\prod_{k=2}^{n}q^{(k-1)(g-1)}\frac{\ \prod_{i=1}^g\bl 1-a_{X,\F_q;k}(i)+q^k\br\ }{(q^k-1)}
\\
=&q^{\binom{n}{2}(g-1)}\prod_{k=1}^{n}\frac{\ \prod_{i=1}^g\bl 1-a_{X,\F_q;k}(i)+q^k\br\ }{(q^k-1)}
\ea
$$
Consequently, note that, by the Riemann hypothesis for $\wh Z_{X,\F_q;k}(T_k)$, 
\be
|a_{X,\F_q;k}(i)|\leq 2\sqrt q^k\qquad(\forall i=1,\ldots g),
\ee 
we arrive at
$$\ba
\prod_{k=1}^{n}\frac{\ \prod_{i=1}^g\bl 1-2\sqrt q^k+q^k\br\ }{(q^k-1)}\leq \b_{X,\F_q;n}(0)
\leq \prod_{k=1}^{n}\frac{\ \prod_{i=1}^g\bl 1+2\sqrt q^k+q^k\br\ }{(q^k-1)}
\ea$$
All these then prove the following
\begin{thm}[Upper and Lower Bounds for $b$-invariants] Let $X$ be an integral  regular projective curve of genus $g$ over a finite field $\F_q$.  Then we have
\be\label{eq23}
\prod_{k=1}^{n}\frac{\ \bl \sqrt q^k-1\br^{2g-1}\ }{(\sqrt q^k+1)}\leq q^{-\binom{n}{2}(g-1)} \b_{X,\F_q;n}(0)
\leq \prod_{k=1}^{n}\frac{\ \bl 1+\sqrt q^k\br^{2g-1}\ }{(\sqrt q^k-1)},
\ee
provided the rank $k$ Riemann hypothesis hold for each $k=2,\ldots,n$.
\end{thm}
\begin{exam} when $X$ is an elliptic curve $E$ over $\F_q$, we have
\be
\prod_{k=1}^{n}\frac{\, \sqrt q^k-1\, }{\sqrt q^k+1}\leq \b_{E,\F_q;n}(0)
\leq \prod_{k=1}^{n}\frac{\, \sqrt q^k+1\, }{\sqrt q^k-1}.
\ee
It is interesting to compare this with Proposition 6 of \cite{WZ1}.
\end{exam}
\begin{exam} In case $n=2$,  we obtain
\be\label{eq23}
q^{g-1}\frac{\  (\sqrt q-1)^{2g-1}(q-1)^{2g-1}\ }{(\sqrt q+1)(q+1)}\leq \b_{X,\F_q;2}(0)
\leq q^{g-1}\frac{\ (\sqrt q+1)^{2g-1}(q+1)^{2g-1}\ }{(\sqrt q-1)(q-1)}.
\ee
Since 
\be
\b_{X,\F_q;2}(0)=q^{g-1}\left(\wh v_2-\frac{1}{q^2-1}\wh v_1^{\,2}\right),
\ee
we obtain
\be
\frac{\  (\sqrt q-1)^{2g-1}(q-1)^{2g-1}\ }{(\sqrt q+1)(q+1)}\leq \wh v_2-\frac{1}{q^2-1}\wh v_1^{\,2}\leq\frac{\  (\sqrt q+1)^{2g-1}(q+1)^{2g-1}\ }{(\sqrt q-1)(q-1)}.
\ee
In this way, we obtain an alternative proof of Lemma\,\ref{cor3.2}.
\end{exam}

\section{Riemann hypothesis for rank three zeta  of a curve over a finite field}

\subsection{Decompose rank three zeta}\label{sec4.1}

Let $X$ be  an integral regular projective curve of genus $g$ over a finite field $\F_q$. By the special uniformity of zetas, we have (up to  a constant factor depending only on $n, g$ and~$q$),
$$
 \ba
 \wh\zeta_{X,\F_q;n}=&\wh\zeta_{X,\F_q}^{\SL_n}(s)\\
 =&\sum_{a=1}^{n}\sum_{\substack{k_1,\ldots,k_p>0\\ k_1+\ldots+k_p=n-a}}\frac{\wh v_{k_1}\ldots\wh v_{k_p}}{\prod_{j=1}^{p-1}(1-q^{k_j+k_{j+1})}} \frac{1}{(1-q^{ns-n+a+k_{p}})}\\
 &\times\wh\zeta_{X,\F_q}(ns-n+a)\sum_{\substack{l_1,\ldots,l_r>0\\ l_1+\ldots+l_r=a-1}}
 \frac{1}{(1-q^{-ns+n-a+1+l_{1}})}\frac{\wh v_{l_1}\ldots\wh v_{l_r}}{\prod_{j=1}^{r-1}(1-q^{l_j+l_{j+1})}}\\
 \ea$$

  In particular, when $n=3$, we have
 $$
 \ba
 \wh\zeta_{X,\F_q;3}(s)=&\wh\zeta_{X,\F_q}^{\SL_3}(s)\\
 =&\left(\frac{\wh v_{1}^{\,2}}{(1-q^{2})(1-q^{3s-1})} 
 +\frac{\wh v_{2}}{(1-q^{3s})}\right)\,\wh\zeta_{X,\F_q}(3s-2)\\
  +&\frac{\wh v_{1}^{\,2}}{(1-q^{3s})(1-q^{-3s+3})}\,\wh\zeta_{X,\F_q}(3s-1) \\
 +&
\left( \frac{\wh v_{2}}{(1-q^{-3s+3})}+
\frac{\wh v_{1}^{\,2}}{(1-q^2)(1-q^{-3s+2})}\right)\, \wh\zeta_{X,\F_q}(3s)\\
 =&\sum_{a=1}^3\wh\zeta_{X/\F_q}^{[a]}(s)
 =:\wh\zeta_{X/\F_q;3}^{\leq2}+\wh\zeta_{X/\F_q;3}^{\geq2}\\
 =&\wh\zeta_{X/\F_q;3}^{\geq2}(s)\left(1+\frac{\wh\zeta_{X/\F_q;3}^{\leq2}(s)}{\wh\zeta_{X/\F_q;3}^{\geq2}(s)}\right)
\ea$$
Here, as before, we have set
 $$\ba
&\wh\zeta_{X/\F_q}^{[1]}(s):=\left(\frac{\wh v_{1}^{\,2}}{(1-q^{2})(1-q^{3s-1})} 
 +\frac{\wh v_{2}}{(1-q^{3s})}\right)\,\wh\zeta_{X,\F_q}(3s-2),\\
&  \wh\zeta_{X/\F_q}^{[2]}(s):=\frac{\wh v_{1}^{\,2}}{(1-q^{3s})(1-q^{-3s+3})}\,\wh\zeta_{X,\F_q}(3s-1), \\
& \wh\zeta_{X/\F_q}^{[3]}(s):=
\left( \frac{\wh v_{2}}{(1-q^{-3s+3})}+
\frac{\wh v_{1}^{\,2}}{(1-q^2)(1-q^{-3s+2})}\right)\, \wh\zeta_{X,\F_q}(3s),\qquad{\rm and}\\
&\wh\zeta_{X/\F_q;3}^{\geq2}(s):=\frac{1}{2}\wh\zeta_{X/\F_q}^{[2]}(s)+\wh\zeta_{X/\F_q}^{[3]}(s)\\
 &=
  \frac{1}{2}\frac{\wh v_{1}^{\,2}}{(1-q^{3s})(1-q^{-3s+3})}\,\wh\zeta_{X,\F_q}(3s-1) 
 +
\left( \frac{\wh v_{2}}{(1-q^{-3s+3})}+
\frac{\wh v_{1}^{\,2}}{(1-q^2)(1-q^{-3s+2})}\right)\, \wh\zeta_{X,\F_q}(3s),\\
 &\wh\zeta_{X/\F_q;3}^{\leq2}(s)
 :=\wh\zeta_{X/\F_q}^{[1]}(s)+\frac{1}{2}\wh\zeta_{X/\F_q}^{[2]}(s)\\
&= \left(\frac{\wh v_{1}^{\,2}}{(1-q^{2})(1-q^{3s-1})} 
 +\frac{1}{2}\frac{\wh v_{2}}{(1-q^{3s})}\right)\,\wh\zeta_{X,\F_q}(3s-2) +\frac{\wh v_{1}^{\,2}}{(1-q^{3s})(1-q^{-3s+3})}\,\wh\zeta_{X,\F_q}(3s-1).
 \ea
 $$
Our strategy, motivated by \cite{S} where a similar result is proved for rank three zeta of the  field $\Q$ of rational, to prove the Riemann Hypothesis for $\wh\zeta_{X/\F_q;3}(s)$ is that first we show the following
\begin{prop}[Riemann Hypothesis for $\wh\zeta_{X/\F_q;3}^{\geq2}(s)$] For an integral  regular projective curve $X$ over $\F_q$, all zeros of $\wh\zeta_{X/\F_q;3}^{\geq2}(s)$ lie on the line $\Re(s)=\frac{1}{3}$.
\end{prop}
Then, based on this proposition, we prove the following
\begin{thm}\label{thm4.2} For an integral  regular projective curve $X$ over $\F_q$, there is no zero of $\wh\zeta_{X,\F_q;3}(s)$ lies in the half plane $\Re(s)<\frac{1}{2}$. 
\end{thm}
Finally, with this theorem, we  prove the following 
\begin{thm}[Riemann Hypothesis in Rank Three]\label{thm4.3} For an integral  regular projective curve $X$ over $\F_q$, all zeros of the rank three non-abelian zeta  $\wh\zeta_{X,\F_q;3}(s)$  and $\SL_3$-zeta $\wh\zeta_{X,\F_q}^{\SL_3}(s)$ of $X/\F_q$ lie on the line $\Re(s)=\frac{1}{2}$.
\end{thm}
\bp
Assuming Theorem\,\ref{thm4.2} whose proof will be given later, our assertion is quite obvious. Indeed,   the functional equation  claims that
\be
\wh\zeta_{X,\F_q;3}(1-s)=\wh\zeta_{X,\F_q;3}(s).
\ee
Hence, there is no zero of $\wh\zeta_{X,\F_q;3}(s)$ lies in the half plane $\Re(s)>\frac{1}{2}$ as well. 
Therefore, all zeros of rank three zeta $\wh\zeta_{X,\F_q;3}(s)$ of $X/\F_q$ lies on the line $\Re(s)=\frac{1}{2}$. 
\ep

\subsection{Estimation on the ratio $\frac{ \wh \zeta_{X,\F_q}(ns-n+a)}{ \wh \zeta_{X,\F_q}(1-ns+n-b)}$ when $a+b=n+1$}
Before proving the results stated above, in this subsection, we will develop some basic techniques to study rank $n$ zetas for curves. We start with the following remarkable 
\begin{prop}\label{prop4.4} Let $X$ be an integral  regular projective curve over $\F_q$. Then for any integers $a$ and $b$ satisfying the condition that $a+b=n+1$, we have
\be
 \left|\frac{ \wh \zeta_{X/\F_q}(ns-n+a)}{ \wh \zeta_{X/\F_q}(1-ns+n-b)}\right|
  \bc
&\hskip -0.3cm>1\qquad{\rm if}\quad |q^{n\s}|<1\\
&\hskip -0.3cm <1\qquad{\rm if}\quad |q^{n\s}|>1.
\ec
 \ee
 where $\s$ is defined by $s=\frac{1}{2}+\s$.
\end{prop}
\bp
The point is to apply the  functional equation for the abelian zeta function
 \be
 \wh\zeta_{X/\F_q}(1-s)= \wh\zeta_{X/\F_q}(s)=\frac{\prod_{i=1}^g(1-\om_{X/\F_q;i}q^{-s})(1-\ov\om_{X/\F_q;i}q^{-s})}{q^{-s(g-1)}(1-q^{-s})(1-q^{1-s})}
 \ee
 to the zeta factor $\wh \zeta_{X,\F_q}(ns-n+b)$ appeared in the denominator of our ratio
 $\dis{\frac{ \wh \zeta_{X,\F_q}(ns-n+a)}{ \wh \zeta_{X,\F_q}(1-ns+n-b)}}$ when $a+b=n+1$.
Accordingly, for $s=\frac{1}{2}+\s$, we have
 $$\ba
 &\frac{\wh \zeta_{X,\F_q}(ns-n+a)\cdot q^{(1-ns+n-b)(g-1)}}{ \wh \zeta_{X,\F_q}(1-ns+n-b)\cdot q^{(ns-n+a)(g-1)}}\\
 &\qquad=\frac{ \wh \zeta_{X,\F_q}(n\s-\frac{n}{2}+a)}{ \wh \zeta_{X,\F_q}(1-n\s+\frac{n}{2}-b)}\cdot q^{(1-2n\s+n-(a+b))(g-1)}\\
 \ea$$
 $$\ba
% =&\left(\prod_{i=1}^g\frac{(1-\om_{X/\F_q,i}q^{-(n\s-\frac{n}{2}+a)})(1-\ov\om_{X/\F_q,i}q^{-(n\s-\frac{n}{2}+a)})}{(1-\om_{X/\F_q,i}q^{-(1-n\s+\frac{n}{2}-b)})(1-\ov\om_{X/\F_q,i}q^{-(1-n\s+\frac{n}{2}-b)})}\right)\frac{(1-q^{-(1-n\s+\frac{n}{2}-b)})(1-q^{1-(1-n\s+\frac{n}{2}-b)})}{(1-q^{-(n\s-\frac{n}{2}+a)})(1-q^{1-(n\s-\frac{n}{2}+a)})}\\
 =&\left(\prod_{i=1}^g\frac{(1-\om_{X/\F_q,i}q^{-n\s+\frac{n}{2}-a})(1-\ov\om_{X/\F_q,i}q^{-n\s+\frac{n}{2}-a})}{(1-\om_{X/\F_q,i}q^{-1+n\s-\frac{n}{2}+b})(1-\ov\om_{X/\F_q,i}q^{-1+n\s-\frac{n}{2}+b})}\right)\frac{(1-q^{-1+n\s-\frac{n}{2}+b})(1-q^{n\s-\frac{n}{2}+b})}{(1-q^{-n\s+\frac{n}{2}-a)})(1-q^{1-n\s+\frac{n}{2}-a})}\\
 %=&q^{g(-2n\s)} \left(\prod_{i=1}^g\frac{(q^{n\s}-\om_{X/\F_q,i}q^{\frac{n}{2}-a})(q^{n\s}-\ov\om_{X/\F_q,i}q^{\frac{n}{2}-a})}{(1-\om_{X/\F_q,i}q^{-1+n\s-\frac{n}{2}+b})(1-\ov\om_{X/\F_q,i}q^{-1+n\s-\frac{n}{2}+b})}\right)\frac{1}{q^{(-2n\s)}}\frac{(1-q^{-1+n\s-\frac{n}{2}+b})(1-q^{n\s-\frac{n}{2}+b)})}{(q^{n\s}-q^{\frac{n}{2}-a)})(q^{n\s}-q^{1+\frac{n}{2}-a)})}\\
 =&q^{-2n\s(g-1)} \left(\prod_{i=1}^g\frac{(q^{n\s}-\om_{X/\F_q,i}q^{\frac{n}{2}-a})(q^{n\s}-\ov\om_{X/\F_q,i}q^{\frac{n}{2}-a})}{(1-\om_{X/\F_q,i}q^{-1+n\s-\frac{n}{2}+b})(1-\ov\om_{X/\F_q,i}q^{-1+n\s-\frac{n}{2}+b})}\right)\\
 &\hskip 4.0cm\times\frac{(1-q^{-1+n\s-\frac{n}{2}+b})(1-q^{n\s-\frac{n}{2}+b})}{(q^{n\s}-q^{\frac{n}{2}-a)})(q^{n\s}-q^{1+\frac{n}{2}-a})}\\
 \ea$$
 To estimate this latest expression, we now introduce a the following key
  \begin{lem}\label{lem3.4} Let $q$ and $\kappa$ be real numbers satisfying $q>1$ and $\kappa\geq 0$. For any complex number $\a,\b$ satisfying $\a\b=q$ and $|\a+\b|\leq q+1$, we have, 
 \be
 \left|w-\a q^\kappa\right|\cdot\left|w-\b q^\kappa\right|=  \left|1-\a q^{\kappa} w\right|\cdot\left|1-\b q^{\kappa} w\right|\cdot \bc>1&\qquad{\rm if}\ |w|<1\\ <1&\qquad{\rm if}\ |w|>1
 \ec
 \ee
 In particular, when $\kappa=0$, we recover Yoshida's inequality in Lemma\,\ref{lem2.1}:
 \be|w-\a|\cdot|w-\b|=|1-\a w|\cdot|1-\b w|\cdot \bc
&\hskip -0.3cm>1\qquad{\rm if}\quad |w|<1\\
&\hskip -0.3cm <1\qquad{\rm if}\quad |w|>1.
\ec
\ee
 \end{lem}

 \bp We start with the following elementary
 \begin{sublem}\label{sublem3.5} For a fixed real number $q>1$, as a function of $x$ in the region $x\geq 0$,
 \be
 f_q(x):=q^{2x+1}+1-q^x(q+1)\geq 0
 \ee
 \end{sublem}
 \bp Clearly, $f(0)=0$. Hence it suffices to verify that
 \be
 f'(x)=\Big(2q^{2x+1}-q^x(q+1)\Big)\log q=q^x\log q\cdot \Big(2q^{x+1}-(q+1)\Big)>0.
 \ee
 But this is a direct consequence of the fact that for $g(x)=2q^{x+1}-(q+1)$,
 \be
 g'(x)=2q^{x+1}\log q>0
\ee
 since $g(0)=q-1>0.$
 \ep
 Back to the proof of the lemma. We follow Yoshida's approach closely. By a direct calculation,
$$\ba
 & |w-\a q^\kappa|^2|w-\b q^\kappa |^2=\Big(w^2-(\a+\b)q^\kappa w+\a\b q^{2\kappa}\Big)\Big(\ov w^2-(\a+\b)q^\kappa \ov w+\a\b q^{2\kappa}\Big)\\
 &=|w|^4-(\a+\b)q^\kappa|w|^2(w+\ov w)+(\a\b)q^{2\kappa}(w^2+\ov w^2)\hskip 5.0cm(a)\\
 &\qquad +(\a+\b)^2 q^{2\kappa}|w|^2-(\a+\b) (\a\b)q^{3\kappa}(w+\ov w)+(\a\b)^2 q^{4\kappa},\qquad{\rm and}\\
 & |1-\a q^\kappa w|^2|1-\b q^\kappa w|^2=\Big(1-(\a+\b)q^\kappa w+\a\b q^{2\kappa}w^2\Big)\Big(1-(\a+\b)q^\kappa \ov w+\a\b q^{2\kappa}\ov w^2\Big)\\
 &=(\a\b)^2q^{4\kappa}|w|^4-(\a+\b)(\a\b)q^{3\kappa}|w|^2(w+\ov w)+(\a\b)q^{2\kappa}(w^2+\ov w^2)\hskip 3.0cm(b)\\
 &\qquad +(\a+\b)^2 q^{2\kappa}|w|^2-(\a+\b) q^{\kappa}(w+\ov w)+1.\\
\ea$$

Subtracting (b) from (a), we get
$$\ba
&|w-\a q^\kappa|^2|w-\b q^\kappa |^2-|1-\a q^\kappa w|^2|1-\b q^\kappa w|^2\\
%&=\Big(1-(\a\b)^2q^{4\kappa}\Big)|w|^4-(\a+\b)q^\kappa\Big(1-(\a\b)q^{2\kappa}\Big)|w|^2(w+\ov w)\\
%&\hskip 1.50cm -(\a+\b) q^\kappa\Big((\a\b)q^{2\kappa}-1\Big)(w+\ov w)+(\a\b)^2 q^{4\kappa}-1\\
%&=\Big(1-(\a\b)q^{2\kappa}\Big)\Big(\Big(1+(\a\b)q^{2\kappa}\Big)|w|^4-(\a+\b)q^\kappa|w|^2(w+\ov w)+(\a+\b) q^{\kappa}(w+\ov w) -(\a\b) q^{2\kappa}-1\Big)\\
%&=\Big(1-(\a\b)q^{2\kappa}\Big)\Big(\Big(1+(\a\b)q^{2\kappa}\Big)\Big(|w|^4-1\Big)-(\a+\b)q^\kappa(w+\ov w)(|w|^2-1)\Big)\Big)\\
&=\Big(1-(\a\b)q^{2\kappa}\Big)\Big(|w|^2-1\Big)\Big(\Big(1+(\a\b)q^{2\kappa}\Big)\Big(|w|^2+1\Big)-(\a+\b)q^\kappa(w+\ov w)\Big)\Big)\\
\ea
$$
Set now $|w|=r$. We claim that
\be\label{eq70}
\Big(1+(\a\b)q^{2\kappa}\Big)\Big(r^2+1\Big)-2r(\a+\b)q^\kappa
=\Big(1+q^{2\kappa+1}\Big)\Big(r^2+1\Big)-2r(\a+\b)q^\kappa>0\qquad(r\not=1)
\ee
since
the discriminant of this degree two polynomial  with real coefficients in $r$ is given by
$$\ba
\De:=&4(\a+\b)^2q^{2\kappa}-4(1+q^{2\kappa+1})^2\leq 4\Big((q+1)^2q^{2\kappa}-\Big(1+2q^{2\kappa+1}+q^{4\kappa+2}\Big)\Big)\\
&=4\Big((q^2+2q+1)q^{2\kappa}-\Big(1+2q^{2\kappa+1}+q^{4\kappa+2}\Big)\Big)=
4\Big((q^2+1)q^{2\kappa}-\Big(1+q^{4\kappa+2}\Big)\Big)\\
&=4f_{q^2}(2\kappa)\leq 0
\ea
$$
by Sublemma\,\ref{sublem3.5}.
Indeed, when $D<0$, the claim is trivial since the leading coefficient $(1+q^{2\kappa+1})$ is strictly positive. In addition, even if this discriminant is zero, which is equivalent to $\kappa=0$ and $\a+\b=q+1$, 
\eqref{eq70} holds as well, since the degree two polynomial becomes $(q+1)(r-1)^2$, which is strictly positive when $r\not=1$.
 \ep
Back to the proof of Proposition\,\ref{prop4.4}.  By the calculation above Lemma\,\ref{lem3.4}, we have
 $$\ba
&\left|\frac{\wh \zeta_{X,\F_q}(ns-n+a)}{ \wh \zeta_{X,\F_q}(1-ns+n-b)}\right|=\left|\frac{ \wh \zeta_{X,\F_q}(n\s-\frac{n}{2}+a)}{ \wh \zeta_{X,\F_q}(1-n\s+\frac{n}{2}-b)}\right|\\
  =&\left|q^{-2n\s(g-2)}\cdot q^{-(1+n-(a+b))(g-1)}\right|\\
  \times&\left(\prod_{i=1}^g\left|\frac{(q^{n\s}-\om_{X/\F_q,i}q^{\frac{n}{2}-a})(q^{n\s}-\ov\om_{X/\F_q,i}q^{\frac{n}{2}-a})}{(1-\om_{X/\F_q,i}q^{-1+n\s-\frac{n}{2}+b})(1-\ov\om_{X/\F_q,i}q^{-1+n\s-\frac{n}{2}+b})}\right|\right)\cdot\left|\frac{(1-q^{-1+n\s-\frac{n}{2}+b})(1-q^{n\s-\frac{n}{2}+b)})}{(q^{n\s}-q^{\frac{n}{2}-a)})(q^{n\s}-q^{1+\frac{n}{2}-a)})}\right|\\
  \ea$$

Now apply Lemma\,\ref{lem3.4} to each pair of factors of the numerator within the product $\prod_{i=1}^g$, namely, with the parameters  $\a=\om_{X/\F_q,i},\ \b=\ov\om_{X/\F_q,i}$ $(i=1,\ldots, g),\ \kappa=\frac{n}{2}~-~a$ and $w=q^{n\s}$, which is applicable thanks to the Riemann hypothesis, or better, the corresponding theorem of Weil, for the Artin zeta function $\wh\zeta_{X/\F_q}(s)$ of the curve $X/\F_q$,  
we have, for $g\geq 2$,\footnote{This is applicable since the genus one rank $n$ Riemann hypothesis has been established in \cite{WZ1}}
  $$\ba
  &\left|\frac{\wh \zeta_{X,\F_q}(ns-n+a)}{ \wh \zeta_{X,\F_q}(1-ns+n-b)}\right|=\left|\frac{ \wh \zeta_{X,\F_q}(n\s-\frac{n}{2}+a)}{ \wh \zeta_{X,\F_q}(1-n\s+\frac{n}{2}-b)}\right|\\
 =&\left| q^{-(1+n-(a+b))(g-1)}\right|\\
 \times &\left(\prod_{i=1}^g\left|\frac{(1-\om_{X/\F_q,i}q^{n\s+\frac{n}{2}-a})(1-\ov\om_{X/\F_q,i}q^{n\s+\frac{n}{2}-a})}{(1-\om_{X/\F_q,i}q^{-1+n\s-\frac{n}{2}+b})(1-\ov\om_{X/\F_q,i}q^{-1+n\s-\frac{n}{2}+b})}\right|\right)\cdot\left|\frac{(1-q^{-1+n\s-\frac{n}{2}+b})(1-q^{n\s-\frac{n}{2}+b)})}{(q^{n\s}-q^{\frac{n}{2}-a)})(q^{n\s}-q^{1+\frac{n}{2}-a)})}\right|\\ 
 &\hskip 1.0cm\times\bc
&\hskip -0.3cm>1\qquad{\rm if}\quad |q^{n\s}|<1\\
&\hskip -0.3cm <1\qquad{\rm if}\quad |q^{n\s}|>1.
\ec
 \ea$$
 In particular, when $a+b=n+1$, we have $\frac{n}{2}-a=-1-\frac{n}{2}+b$  and hence the denominator and  numerators are identical for each pair of factors within the product $\prod_{i=1}^g$. Consequently,
  $$\ba
 &\left|\frac{ \wh \zeta_{X,\F_q}(ns-n+a)}{ \wh \zeta_{X,\F_q}(1-ns+n-b)}\right|
  =\left|\frac{(1-q^{n\s+\frac{n}{2}-a})(1-q^{n\s+\frac{n}{2}-a+1)})}{(q^{n\s}-q^{\frac{n}{2}-a)})(q^{n\s}-q^{1+\frac{n}{2}-a)})}\right|\cdot\bc
&\hskip -0.3cm>1\qquad{\rm if}\quad |q^{n\s}|<1\\
&\hskip -0.3cm <1\qquad{\rm if}\quad |q^{n\s}|>1.
\ec
 \ea$$
Now,  by applying Lemma\,\ref{lem3.4} again to the numerators involved, with parameters $\a=~1$, $b=q$, $\kappa=\frac{n}{2}-a$ and $w=q^{n\s}$,
we get
  \be
 \left|\frac{ \wh \zeta_{X,\F_q}(ns-n+a)}{ \wh \zeta_{X,\F_q}(1-ns+n-b)}\right|
  \bc
&\hskip -0.3cm>1\qquad{\rm if}\quad |q^{n\s}|<1\\
&\hskip -0.3cm <1\qquad{\rm if}\quad |q^{n\s}|>1.
\ec\qquad (a+b=n+1,\ s=:\frac{1}{2}+\s)
 \ee
 as wanted.
 \ep

\subsection{Estimate the ratio $\frac{ \wh\zeta_{X,\F_q;n}^{[a]}(s)}{ \wh\zeta_{X,\F_q;n}^{[b]}(s)}$ when $b=a+1$}

By definition,
  $$\ba
\frac{ \wh\zeta_{X,\F_q;n}^{[a]}(s)}{ \wh\zeta_{X,\F_q;n}^{[b]}(s)}
 =&\frac{\sum_{\substack{k_1,\ldots,k_p>0\\ k_1+\ldots+k_p=n-a}}\frac{\wh v_{k_1}\ldots\wh v_{k_p}}{\prod_{j=1}^{p-1}(1-q^{k_j+k_{j+1})}} \frac{1}{(1-q^{-n+a+k_{p}}/T)}
\sum_{\substack{l_1,\ldots,l_r>0\\ l_1+\ldots+l_r=a-1}}
 \frac{1}{(1-Tq^{n-a+1+l_{1}})}\frac{\wh v_{l_1}\ldots\wh v_{l_r}}{\prod_{j=1}^{p-1}(1-q^{l_j+l_{j+1})}}
}{\sum_{\substack{k_1,\ldots,k_p>0\\ k_1+\ldots+k_p=n-b}}\frac{\wh v_{k_1}\ldots\wh v_{k_p}}{\prod_{j=1}^{p-1}(1-q^{k_j+k_{j+1})}} \frac{1}{(1-q^{-n+b+k_{p}}/T)}
\sum_{\substack{l_1,\ldots,l_r>0\\ l_1+\ldots+l_r=b-1}}
 \frac{1}{(1-Tq^{n-b+1+l_{1}})}\frac{\wh v_{l_1}\ldots\wh v_{l_r}}{\prod_{j=1}^{p-1}(1-q^{l_j+l_{j+1})}}
}\\
&\times\frac{ \wh \zeta_{X,\F_q}(ns-n+a)}{ \wh \zeta_{X,\F_q}(ns-n-b)}\\
 \ea$$
 Write $s=\frac{1}{2}+\s$, then the zeta factor becomes
 $$\ba
 &\frac{ \wh \zeta_{X,\F_q}(ns-n+a)}{ \wh \zeta_{X,\F_q}(ns-n+b)}=\frac{ \wh \zeta_{X,\F_q}(n\s-\frac{n}{2}+a)}{ \wh \zeta_{X,\F_q}(n\s-\frac{n}{2}+b)}\\
 =&\frac{q^{-(n\s-\frac{n}{2}+b)(g-1)}}{q^{-(n\s-\frac{n}{2}+a)(g-1)}}\\
 &\times\left(\prod_{i=1}^g\frac{(1-\om_{X/\F_q,i}q^{-(n\s-\frac{n}{2}+a)})(1-\ov\om_{X/\F_q,i}q^{-(n\s-\frac{n}{2}+a)})}{(1-\om_{X/\F_q,i}q^{-(n\s-\frac{n}{2}+b)})(1-\ov\om_{X/\F_q,i}q^{-(n\s-\frac{n}{2}+b)})}\right)\cdot \frac{(1-q^{-(n\s-\frac{n}{2}+b)})(1-q^{1-(n\s-\frac{n}{2}+b)})}{(-q^{-(n\s-\frac{n}{2}+a)})(1-q^{1-(n\s-\frac{n}{2}+a)})}\\
% =&\left(\prod_{i=1}^g\frac{(1-\om_{X/\F_q,i}q^{-n\s+\frac{n}{2}-a})(1-\ov\om_{X/\F_q,i}q^{-n\s+\frac{n}{2}-a})}{(1-\om_{X/\F_q,i}q^{-n\s+\frac{n}{2}-b})(1-\ov\om_{X/\F_q,i}q^{-n\s+\frac{n}{2}-b})}\right)\cdot \frac{(1-q^{-n\s+\frac{n}{2}-b})(1-q^{1-n\s+\frac{n}{2}-b)})}{(1-q^{-n\s+\frac{n}{2}-a)})(1-q^{1-n\s+\frac{n}{2}-a)})}\cdot q^{(a-b)(g-1)}\\
 %  =&\left(\prod_{i=1}^g\frac{(1-\om_{X/\F_q,i}q^{-n\s+\frac{n}{2}-a})(1-\ov\om_{X/\F_q,i}q^{-n\s+\frac{n}{2}-a})}{(1-\om_{X/\F_q,i}q^{-n\s+\frac{n}{2}-a-1})(1-\ov\om_{X/\F_q,i}q^{-n\s+\frac{n}{2}-a-1})}\right)\cdot\frac{(1-q^{-n\s+\frac{n}{2}-a-1})(1-q^{-n\s+\frac{n}{2}-a})}{(1-q^{-n\s+\frac{n}{2}-a})(1-q^{1-n\s+\frac{n}{2}-a)})}\cdot q^{-(g-1)}\\
    %&\hskip 10.0cm (b=a+1)\\
  %  =&\left(\prod_{i=1}^g\frac{(1-\om_{X/\F_q,i}q^{-n\s+\frac{n}{2}-a})(1-\ov\om_{X/\F_q,i}q^{-n\s+\frac{n}{2}-a})}{(1-\frac{1}{\ov\om_{X/\F_q,i}}q^{-n\s+\frac{n}{2}-a})(1-\frac{1}{\om_{X/\F_q,i}}q^{-n\s+\frac{n}{2}-a})}\right)\cdot \frac{(1-q^{-n\s+\frac{n}{2}-a-1})}{(1-q^{1-n\s+\frac{n}{2}-a)})}\cdot q^{-(g-1)}\\
 % &\hskip 7.0cm (\om_{X/\F_q,i}\ov\om_{X,\F_q}=q\qquad i=1,\ldots,g)\\
   =&\left(\prod_{i=1}^g\frac{(1-\om_{X/\F_q,i}q^{-n\s+\frac{n}{2}-a})(1-\ov\om_{X/\F_q,i}q^{-n\s+\frac{n}{2}-a})}{({\ov\om_{X/\F_q,i}}-q^{-n\s+\frac{n}{2}-a})({\om_{X/\F_q,i}}-q^{-n\s+\frac{n}{2}-a})}\right)\cdot \frac{(q-q^{-n\s+\frac{n}{2}-a})}{(1-q^{1-n\s+\frac{n}{2}-a)})}\\
  &\hskip 3.40cm ({\rm since}\ b=a+1\qan \om_{X/\F_q,i}\ov\om_{X,\F_q}=q\quad\forall  i=1,\ldots,g)\\
 \ea$$
Therefore, by applying Yoshida's lemma to the factors in the denominator with the parameter
$\a=\om_{X/\F_q,i}$,\ $\,\b=\ov \om_{X/\F_q,i}$ and $w=q^{-n\s+\frac{n}{2}-a}$, we have
\be
\left|\frac{ \wh \zeta_{X,\F_q}(n\s-\frac{n}{2}+a)}{ \wh \zeta_{X,\F_q}(n\s-\frac{n}{2}+b)}\right|
%=&\left|\prod_{i=1}^g\frac{(1-\om_{X/\F_q,i}q^{-n\s+\frac{n}{2}-a})(1-\ov\om_{X/\F_q,i}q^{-n\s+\frac{n}{2}-a})}{({\ov\om_{X/\F_q,i}}-q^{-n\s+\frac{n}{2}-a})({\om_{X/\F_q,i}}-q^{-n\s+\frac{n}{2}-a})}\right|\cdot \left|\frac{(q-q^{-n\s+\frac{n}{2}-a})}{(1-q^{1-n\s+\frac{n}{2}-a)})}\right|\\
=\left|\frac{(q-q^{-n\s+\frac{n}{2}-a})}{(1-q^{1-n\s+\frac{n}{2}-a)})}\right|\cdot\bc
>1&\qquad\left|q^{-n\s+\frac{n}{2}-a}\right|>1\\
<1&\qquad\left|q^{-n\s+\frac{n}{2}-a}\right|<1
\ec\quad(b=a+1)
\ee
That is to say, we have proved the following
\begin{prop}\label{prop4.7} For an integral  regular projective curve  $X$ over $\F_q$, we have
\be
\left|\frac{ \wh \zeta_{X,\F_q}(n\s-\frac{n}{2}+a)}{ \wh \zeta_{X,\F_q}(n\s-\frac{n}{2}+b)}\right|\\
= \left|\frac{(q-q^{-n\s+\frac{n}{2}-a})}{(1-q^{1-n\s+\frac{n}{2}-a)})}\right|\cdot\bc
>1&\qquad \left|q^{n\s}\right|<\left|q^{\frac{n}{2}-a}\right|\\
<1&\qquad \left|q^{n\s}\right|>\left|q^{\frac{n}{2}-a}\right|
\ec\qquad(b=a+1)
\ee 
\end{prop}
This implies that, when $b=a+1$,
$$\ba
&\left|\frac{ \wh\zeta_{X,\F_q;n}^{[a]}(s)}{ \wh\zeta_{X,\F_q;n}^{[b]}(s)}\right|\\
 =&\left|\frac{\sum_{\substack{k_1,\ldots,k_p>0\\ k_1+\ldots+k_p=n-a}}\frac{\wh v_{k_1}\ldots\wh v_{k_p}}{\prod_{j=1}^{p-1}(1-q^{k_j+k_{j+1}})} \frac{1}{(1-q^{-n+a+k_{p}+\frac{n}{2}+n\s})}
\sum_{\substack{l_1,\ldots,l_r>0\\ l_1+\ldots+l_r=a-1}}
 \frac{1}{(1-q^{-\frac{n}{2}-n\s+n-a+1+l_{1}})}\frac{\wh v_{l_1}\ldots\wh v_{l_r}}{\prod_{j=1}^{p-1}(1-q^{l_j+l_{j+1}})}
}{\sum_{\substack{k_1,\ldots,k_p>0\\ k_1+\ldots+k_p=n-a-1}}\frac{\wh v_{k_1}\ldots\wh v_{k_p}}{\prod_{j=1}^{p-1}(1-q^{k_j+k_{j+1}})} \frac{1}{(1-q^{-n+a+1+k_{p}+\frac{n}{2}+n\s})}
\sum_{\substack{l_1,\ldots,l_r>0\\ l_1+\ldots+l_r=a}}
 \frac{1}{(1-q^{-\frac{n}{2}-n\s+n-a+l_{1}})}\frac{\wh v_{l_1}\ldots\wh v_{l_r}}{\prod_{j=1}^{p-1}(1-q^{l_j+l_{j+1}})}
}\right|\\
&\times \left|\frac{(q-q^{-n\s+\frac{n}{2}-a})}{(1-q^{1-n\s+\frac{n}{2}-a)})}\right|\cdot\bc
>1&\qquad \left|q^{n\s}\right|<\left|q^{\frac{n}{2}-a}\right|\\
<1&\qquad \left|q^{n\s}\right|>\left|q^{\frac{n}{2}-a}\right|
\ec\\
=&\left|\frac{\sum_{\substack{k_1,\ldots,k_p>0\\ k_1+\ldots+k_p=n-a}}\frac{\wh v_{k_1}\ldots\wh v_{k_p}}{\prod_{j=1}^{p-1}(1-q^{k_j+k_{j+1}})} \frac{1}{(1-q^{(n\s-\frac{n}{2}+a)+k_{p}})}
\sum_{\substack{l_1,\ldots,l_r>0\\ l_1+\ldots+l_r=a-1}}
 \frac{1}{(1-q^{-(n\s-\frac{n}{2}+a)+1+l_{1}})}\frac{\wh v_{l_1}\ldots\wh v_{l_r}}{\prod_{j=1}^{p-1}(1-q^{l_j+l_{j+1}})}
}{\sum_{\substack{k_1,\ldots,k_p>0\\ k_1+\ldots+k_p=n-a-1}}\frac{\wh v_{k_1}\ldots\wh v_{k_p}}{\prod_{j=1}^{p-1}(1-q^{k_j+k_{j+1}})} \frac{1}{(1-q^{(n\s-\frac{n}{2}+a)+1+k_{p}})}
\sum_{\substack{l_1,\ldots,l_r>0\\ l_1+\ldots+l_r=a}}
 \frac{1}{(1-q^{-(n\s-\frac{n}{2}+a)+l_{1}})}\frac{\wh v_{l_1}\ldots\wh v_{l_r}}{\prod_{j=1}^{p-1}(1-q^{l_j+l_{j+1}})}
}\right|\\
&\times  \left|\frac{(q-q^{-(n\s-\frac{n}{2}+a)})}{(1-q^{1-(n\s-\frac{n}{2}+a)})}\right|\cdot\bc
>1&\qquad \left|q^{n\s}\right|<\left|q^{\frac{n}{2}-a}\right|\\
<1&\qquad \left|q^{n\s}\right|>\left|q^{\frac{n}{2}-a}\right|
\ec\\
\ea$$
So we are lead to consider the norm of the following rational function in the first big factor
$$\ba
&r_{X/\F_q;n,a}(\s)\\
:=&\frac{\sum_{\substack{k_1,\ldots,k_p>0\\ k_1+\ldots+k_p=n-a}}\frac{\wh v_{k_1}\ldots\wh v_{k_p}}{\prod_{j=1}^{p-1}(1-q^{k_j+k_{j+1}})} \frac{1}{(1-q^{(n\s-\frac{n}{2}+a)+k_{p}})}
\sum_{\substack{l_1,\ldots,l_r>0\\ l_1+\ldots+l_r=a-1}}
 \frac{1}{(1-q^{-(n\s-\frac{n}{2}+a)+1+l_{1}})}\frac{\wh v_{l_1}\ldots\wh v_{l_r}}{\prod_{j=1}^{p-1}(1-q^{l_j+l_{j+1}})}
}{\sum_{\substack{k_1,\ldots,k_p>0\\ k_1+\ldots+k_p=n-a-1}}\frac{\wh v_{k_1}\ldots\wh v_{k_p}}{\prod_{j=1}^{p-1}(1-q^{k_j+k_{j+1}})} \frac{1}{(1-q^{(n\s-\frac{n}{2}+a)+1+k_{p}})}
\sum_{\substack{l_1,\ldots,l_r>0\\ l_1+\ldots+l_r=a}}
 \frac{1}{(1-q^{-(n\s-\frac{n}{2}+a)+l_{1}})}\frac{\wh v_{l_1}\ldots\wh v_{l_r}}{\prod_{j=1}^{p-1}(1-q^{l_j+l_{j+1}})}
}\\
=:&\frac{f_{X/\F_q;n,a}(\s)\cdot g_{X/\F_q;n,a}(\s)}{f_{X/\F_q;n,a+1}(\s)\cdot g_{X/\F_q;n,a+1}(\s)},\ea$$
where in the last step, we have set accordingly
\be\label{eq73}
f_{X/\F_q;n,a}(\s):=\sum_{\substack{k_1,\ldots,k_p>0\\ k_1+\ldots+k_p=n-a}}\frac{\wh v_{k_1}\ldots\wh v_{k_p}}{\prod_{j=1}^{p-1}(1-q^{k_j+k_{j+1}})} \frac{1}{(1-q^{(n\s-\frac{n}{2}+a)+k_{p}})}
\ee
and 
\be\label{eq74}
g_{X/\F_q;n,a}(\s):=\sum_{\substack{l_1,\ldots,l_r>0\\ l_1+\ldots+l_r=a-1}}
 \frac{1}{(1-q^{-(n\s-\frac{n}{2}+a)+1+l_{1}})}\frac{\wh v_{l_1}\ldots\wh v_{l_r}}{\prod_{j=1}^{p-1}(1-q^{l_j+l_{j+1}})}
\ee
Therefore,
\be\label{eq75}
\left|\frac{ \wh\zeta_{X,\F_q;n}^{[a]}(s)}{ \wh\zeta_{X,\F_q;n}^{[b]}(s)}\right|
=\frac{f_{X/\F_q;n,a}(\s)\cdot g_{X/\F_q;n,a}(\s)}{f_{X/\F_q;n,a+1}(\s)\cdot g_{X/\F_q;n,a+1}(\s)}
 \left|\frac{(q-q^{-(n\s-\frac{n}{2}+a)})}{(1-q^{1-(n\s-\frac{n}{2}+a)})}\right|
 \cdot\bc
>1&\qquad \left|q^{n\s}\right|<\left|q^{\frac{n}{2}-a}\right|\\
<1&\qquad \left|q^{n\s}\right|>\left|q^{\frac{n}{2}-a}\right|\ec
\ee

We end this subsection with the following comments. It appears to be very tempting  to apply Yoshida's lemma to the factor in the middle, namely, $ \left|\frac{(q-q^{-(n\s-\frac{n}{2}+a)})}{(1-q^{1-(n\s-\frac{n}{2}+a)})}\right|$ with parameters $\a=1,\ \b=q$ and $w=q^{-(n\s-\frac{n}{2}+a)}$. Unfortunately, this would only result  inequalities in opposite directions. Nevertheless, as to be seen in the subsection below, there will be  a nice total cancelation on this middle factor from the factors in the first group on the ratios of $f$ and $g$'s.

\subsection{The Riemann hypothesis for $ \wh\zeta_{X,\F_q;3}^{\,\geq 2}(s)$}
In this subsection, we prove the following
\begin{prop}[Riemann Hypothesis for $\wh\zeta_{X,\F_q;3}^{\,\geq 2}(s)$] Let $X$ be  an integral  regular projective curve over $\F_q$. Then all zeros of $ \wh\zeta_{X,\F_q;3}^{\,\geq 2}(s)$ lie on the line $\Re(s)=\frac{1}{3}$.
\end{prop}
\bp
Indeed, from the definitions in the previous subsection, particularly, in \eqref{eq73} and  \eqref{eq74} taking parameters $n=3$, $a=2$ and $a=3$, we have
$$\ba
f_{X/\F_q;3,2}(\s):=&\sum_{\substack{k_1,\ldots,k_p>0\\ k_1+\ldots+k_p=1}}\frac{\wh v_{k_1}\ldots\wh v_{k_p}}{\prod_{j=1}^{p-1}(1-q^{k_j+k_{j+1}})} \frac{1}{(1-q^{(3\s-\frac{3}{2}+2)+k_{p}})}\\
=&\sum_{p=1,\, k_p=1}\wh v_{k_p}\frac{1}{(1-q^{(3\s-\frac{3}{2}+2)+k_{p}})}=\wh v_1\frac{1}{(1-q^{(3\s-\frac{3}{2}+2)+1})}=\frac{\wh v_1}{1-q^{3\s+\frac{3}{2}}}
\ea
$$
and 
$$
f_{X/\F_q;3,2+1}(\s):=\sum_{\substack{k_1,\ldots,k_p>0\\ k_1+\ldots+k_p=3-3}}\frac{\wh v_{k_1}\ldots\wh v_{k_p}}{\prod_{j=1}^{p-1}(1-q^{k_j+k_{j+1}})} \frac{1}{(1-q^{(3\s-\frac{3}{2}+3)+k_{p}})}=1
$$
Similarly,
$$\ba
g_{X/\F_q;3,2}(\s):=&\sum_{\substack{l_1,\ldots,l_r>0\\ l_1+\ldots+l_r=2-1}}
 \frac{1}{(1-q^{-(3\s-\frac{3}{2}+2)+1+l_{1}})}\frac{\wh v_{l_1}\ldots\wh v_{l_r}}{\prod_{j=1}^{p-1}(1-q^{l_j+l_{j+1}})}\\
 =&\sum_{r=1,l_r=1} \frac{1}{(1-q^{-(3\s-\frac{3}{2}+2)+1+l_{1}})}\wh v_{l_1}
 =
 \frac{\wh v_1}{1-q^{-3\s+\frac{3}{2}}}
\ea$$
and 
$$\ba
g_{X/\F_q;3,2+1}(\s)
:=&\sum_{\substack{l_1,\ldots,l_r>0\\ l_1+\ldots+l_r=3-1}}
 \frac{1}{(1-q^{-(3\s-\frac{3}{2}+3)+1+l_{1}})}\frac{\wh v_{l_1}\ldots\wh v_{l_r}}{\prod_{j=1}^{p-1}(1-q^{l_j+l_{j+1}})}\\
 %=&\sum_{\substack{l_1,\ldots,l_r>0\\ l_1+\ldots+l_r=2}}
% \frac{1}{(1-q^{-(3\s-\frac{3}{2}+3)+1+l_{1}})}\frac{\wh v_{l_1}\ldots\wh v_{l_r}}{\prod_{j=1}^{p-1}(1-q^{l_j+l_{j+1}})}\\
 =&\sum_{r=1,l_r=2} \frac{1}{(1-q^{-(3\s-\frac{3}{2}+3)+1+l_{1}})}\wh v_{l_1}+\sum_{\substack{r=2\\ l_1= l_r=1}} \frac{1}{(1-q^{-(3\s-\frac{3}{2}+3)+1+l_{1}})}\frac{\wh v_{l_1}\wh v_{l_r}}{\prod_{j=1}^{2-1}(1-q^{l_j+l_{j+1}})}\\
 %=&\frac{\wh v_{2}}{(1-q^{-(3\s-\frac{3}{2}+3)+1+2})}+\frac{1}{(1-q^{-(3\s-\frac{3}{2}+3)+1+1})}%\frac{\wh v_{1}^2}{1-q^2}\\
 =&\frac{\wh v_{2}}{(1-q^{-3\s+\frac{3}{2}})}+\frac{1}{(1-q^{-3\s+\frac{3}{2}-1})}\frac{\wh v_{1}^2}{1-q^2}\\
 \ea
 $$
 Therefore,
 $$\ba
 r_{X/\F_q;3,2}(\s)=&\frac{f_{X/\F_q;3,2}(\s)\cdot g_{X/\F_q;3,2}(\s)}{f_{X/\F_q;3,2+1}(\s)\cdot g_{X/\F_q;3,2+1}(\s)}\\
 %=&\frac{\frac{\wh v_1}{1-q^{3\s+\frac{3}{2}}}
%\cdot \frac{\wh v_1}{1-q^{-3\s+\frac{3}{2}}} }{1\cdot \left(\frac{\wh v_{2}}{(1-q^{-3\s+\frac{3}{2}})}+\frac{1}{(1-q^{-3\s+\frac{3}{2}-1})}\frac{\wh v_{1}^2}{1-q^2}\right)}\\
%=&\frac{\wh v_1^{\,2}}
%{\left({(1-q^{-3\s+\frac{3}{2}-1})\wh v_{2}}+(1-q^{-3\s+\frac{3}{2}})\frac{\wh v_{1}^2}{1-q^2}\right)}\cdot \frac{(1-q^{-3\s+\frac{3}{2}})(1-q^{-3\s+\frac{3}{2}-1})}{(1-q^{3\s+\frac{3}{2}})(1-q^{-3\s+\frac{3}{2}})} \\
%=&\frac{\wh v_1^{\,2}}
%{\left(\Big({\wh v_{2}}+\frac{\wh v_{1}^2}{1-q^2}\Big)-q^{-3\s+\frac{3}{2}-1}\Big(\wh v_{2}+\frac{q\,\wh v_{1}^2}{1-q^2}\Big)\right)}\cdot \frac{\cancel{(1-q^{-3\s+\frac{3}{2}})}(1-q^{-3\s+\frac{3}{2}-1})}{(1-q^{3\s+\frac{3}{2}})\cancel{(1-q^{-3\s+\frac{3}{2}})}} \\
=&\frac{1}
{\left(\Big(\frac{\wh v_{2}}{\wh v_1^{\,2}}+\frac{1}{1-q^2}\Big)-q^{-3\s+\frac{3}{2}-1}\Big(\frac{\wh v_{2}}{\wh v_1^{\,2}}+\frac{q}{1-q^2}\Big)\right)}\cdot \frac{1-q^{-3\s+\frac{3}{2}-1}}{1-q^{3\s+\frac{3}{2}}} 
 \ea$$
Consequently, by the result in the previous subsection, particularly, by \eqref{eq75}, we have
$$\ba
&\left|\frac{ \wh\zeta_{X,\F_q;3}^{[2]}(s)}{ \wh\zeta_{X,\F_q;3}^{[3]}(s)}\right|
 =\Big|r_{X/\F_q;3,2}(\s)\Big| \cdot  \left|\frac{(q-q^{-(3\s-\frac{3}{2}+2)})}{(1-q^{1-(3\s-\frac{3}{2}+2)})}\right|\cdot\bc
>1&\qquad \left|q^{3\s}\right|<\left|q^{\frac{3}{2}-2}\right|\\
<1&\qquad \left|q^{3\s}\right|>\left|q^{\frac{3}{2}-2}\right|
\ec\\
%=&\frac{1}
%{\left|\Big(\frac{\wh v_{2}}{\wh v_1^{\,2}}+\frac{1}{1-q^2}\Big)-q^{-3\s+\frac{3}{2}-1}\Big(\frac{\wh v_{2}}{\wh v_1^{\,2}}+\frac{q}{1-q^2}\Big)\right|}\cdot \left|\frac{\xcancel{(1-q^{-3\s+\frac{3}{2}-1})}}{(1-q^{3\s+\frac{3}{2}})}\right| \cdot \left|\frac{(q-q^{-(3\s-\frac{3}{2}+2)})}{\xcancel{(1-q^{1-(3\s-\frac{3}{2}+2)})}}\right|\\
%&\times\bc
%>1&\qquad \left|q^{3\s}\right|<\left|q^{-\frac{1}{2}}\right|\\
%<1&\qquad \left|q^{3\s}\right|>\left|q^{-\frac{1}{2}}\right|
%\ec\\
%=&\frac{q}
%{\left|\Big(\frac{\wh v_{2}}{\wh v_1^{\,2}}+\frac{1}{1-q^2}\Big)-q^{-3\s+\frac{3}{2}-1}\Big(\frac{\wh v_{2}}{\wh v_1^{\,2}}+\frac{q}{1-q^2}\Big)\right|}\cdot \left|\frac{1}{(1-q^{3\s+\frac{3}{2}})}\right| \cdot \Big|1-q^{-3\s-\frac{3}{2}}\Big|\\
%&\times\bc
%>1&\qquad \left|q^{3\s}\right|<\left|q^{-\frac{1}{2}}\right|\\
%<1&\qquad \left|q^{3\s}\right|>\left|q^{-\frac{1}{2}}\right|
%\ec\\
%=&\frac{q\Big|q^{-3\s-\frac{3}{2}}\Big|}
%{\left|\Big(\frac{\wh v_{2}}{\wh v_1^{\,2}}+\frac{1}{1-q^2}\Big)-q^{-3\s+\frac{3}{2}-1}\Big(\frac{\wh v_{2}}{\wh v_1^{\,2}}+\frac{q}{1-q^2}\Big)\right|}\cdot \left|\frac{1}{\cancel{(1-q^{3\s+\frac{3}{2}})}}\right| \cdot \Big|\cancel{1-q^{3\s+\frac{3}{2}}}\Big|\\
%&\times\bc
%>1&\qquad \left|q^{3\s}\right|<\left|q^{-\frac{1}{2}}\right|\\
%<1&\qquad \left|q^{3\s}\right|>\left|q^{-\frac{1}{2}}\right|
%\ec\\
=&\frac{\Big|q^{-3\s-\frac{1}{2}}\Big|}
{\left|\Big(\frac{\wh v_{2}}{\wh v_1^{\,2}}+\frac{1}{1-q^2}\Big)-q^{-3\s+\frac{3}{2}-1}\Big(\frac{\wh v_{2}}{\wh v_1^{\,2}}+\frac{q}{1-q^2}\Big)\right|}\cdot \bc
>1&\qquad \left|q^{3\s}\right|<\left|q^{-\frac{1}{2}}\right|\\
<1&\qquad \left|q^{3\s}\right|>\left|q^{-\frac{1}{2}}\right|
\ec\\
=&\frac{1}
{\left|\Big(\frac{\wh v_{2}}{\wh v_1^{\,2}}+\frac{1}{1-q^2}\Big)-q^{-3\s+\frac{3}{2}-1}\Big(\frac{\wh v_{2}}{\wh v_1^{\,2}}+\frac{q}{1-q^2}\Big)\right|}\cdot \bc
>1&\qquad \left|q^{3\s}\right|<\left|q^{-\frac{1}{2}}\right|\\
<1&\qquad \left|q^{3\s}\right|>\left|q^{-\frac{1}{2}}\right|
\ec\\
\ea$$

\begin{lem}\label{lem4.8}  we have
\be
\frac{1/2}
{\left|\Big(\frac{\wh v_{2}}{\wh v_1^{\,2}}+\frac{1}{1-q^2}\Big)-q^{-3\s+\frac{3}{2}-1}\Big(\frac{\wh v_{2}}{\wh v_1^{\,2}}+\frac{q}{1-q^2}\Big)\right|} \bc
>1&\qquad \left|q^{3\s}\right|<\left|q^{-\frac{1}{2}}\right|\\
<1&\qquad \left|q^{3\s}\right|>\left|q^{-\frac{1}{2}}\right|
\ec
\ee
\end{lem}

\bp
To see it clearly, set $\s=-\frac{1}{6}+\tau$ so that $3\s=-\frac{1}{2}+3\tau$, and let
$$\ba
D_{X/\F_q,3,2}(\tau):=&2\left|\Big(\frac{\wh v_{2}}{\wh v_1^{\,2}}+\frac{1}{1-q^2}\Big)-q^{-3\s+\frac{3}{2}-1}\Big(\frac{\wh v_{2}}{\wh v_1^{\,2}}+\frac{q}{1-q^2}\Big)\right|\\
=&2\left|\Big(\frac{\wh v_{2}}{\wh v_1^{\,2}}+\frac{1}{1-q^2}\Big)-q^{-3\tau}\Big(q\frac{\wh v_{2}}{\wh v_1^{\,2}}+\frac{q^2}{1-q^2}\Big)\right|\\
=&2q\cdot\Big(\frac{\wh v_{2}}{\wh v_1^{\,2}}+\frac{q}{1-q^2}\Big)\cdot 
{\left|\frac{\Big(\frac{\wh v_{2}}{\wh v_1^{\,2}}+\frac{1}{1-q^2}\Big)}{\Big(q\frac{\wh v_{2}}{\wh v_1^{\,2}}+\frac{q^2}{1-q^2}\Big)}-q^{-3\tau}\right|}
\ea
$$
where in the last step, we have used Corollary\,\ref{cor3.2} that 
\be
\frac{\wh v_{2}}{\wh v_1^{\,2}}+\frac{q}{1-q^2}>0.
\ee
Note that the condition $\left|q^{3\s}\right|<\left|q^{-\frac{1}{2}}\right|
$ is equivalent to $\left|q^{3\tau}\right|<1$, and similarly for the opposite direction.
Hence it suffices to verify that  
\be
2q\cdot\Big(\frac{\wh v_{2}}{\wh v_1^{\,2}}+\frac{q}{1-q^2}\Big)\cdot 
{\left|\frac{\Big(\frac{\wh v_{2}}{\wh v_1^{\,2}}+\frac{1}{1-q^2}\Big)}{\Big(q\frac{\wh v_{2}}{\wh v_1^{\,2}}+\frac{q^2}{1-q^2}\Big)}-q^{-3\tau}\right|}=
D_{X/\F_q,3,2}(\tau)\bc
<1&\qquad \left|q^{3\tau}\right|<1\\
>1&\qquad \left|q^{3\tau}\right|>1
\ec
\ee
In other words, for $w=q^{3\tau}$, we have to show that
\begin{enumerate}
\item [(1)] When $|w|<1$, then
$w$ should be contained inside the disc of radius $\dis{\frac{1}
{2q\Big(\frac{\wh v_{2}}{\wh v_1^{\,2}}+\frac{q}{1-q^2}\Big)}}$ centered at $\dis{\frac{\frac{\wh v_{2}}{\wh v_1^{\,2}}+\frac{1}{1-q^2}}{q\frac{\wh v_{2}}{\wh v_1^{\,2}}+\frac{q^2}{1-q^2}}}$; and 
\item [(2)] When $|w|>1$, then $w$ should be totally located outside the disc of radius $\dis{\frac{1}
{2q\Big(\frac{\wh v_{2}}{\wh v_1^{\,2}}+\frac{q}{1-q^2}\Big)}}$ centered at $\dis{\frac{\frac{\wh v_{2}}{\wh v_1^{\,2}}+\frac{1}{1-q^2}}{q\frac{\wh v_{2}}{\wh v_1^{\,2}}+\frac{q^2}{1-q^2}}}$.
\end{enumerate}
An elementary discussion implies that this would happen  if the disc of radius $\dis{\frac{1}
{2q\Big(\frac{\wh v_{2}}{\wh v_1^{\,2}}+\frac{q}{1-q^2}\Big)}}$ centered at $\dis{\frac{\frac{\wh v_{2}}{\wh v_1^{\,2}}+\frac{1}{1-q^2}}{q\frac{\wh v_{2}}{\wh v_1^{\,2}}+\frac{q^2}{1-q^2}}}$ is totally contained in the unit disc $|w|<1$. 
Since $\dis{\frac{\frac{\wh v_{2}}{\wh v_1^{\,2}}+\frac{1}{1-q^2}}{q\frac{\wh v_{2}}{\wh v_1^{\,2}}+\frac{q^2}{1-q^2}}}\in \R$ is a real number, this means that
\be
\bc -\dis{\frac{1}
{2q\Big(\frac{\wh v_{2}}{\wh v_1^{\,2}}+\frac{q}{1-q^2}\Big)}}+\dis{\frac{\frac{\wh v_{2}}{\wh v_1^{\,2}}+\frac{1}{1-q^2}}{q\frac{\wh v_{2}}{\wh v_1^{\,2}}+\frac{q^2}{1-q^2}}}>-1\\
\dis{\frac{1}
{2q\Big(\frac{\wh v_{2}}{\wh v_1^{\,2}}+\frac{q}{1-q^2}\Big)}}+\dis{\frac{\frac{\wh v_{2}}{\wh v_1^{\,2}}+\frac{1}{1-q^2}}{q\frac{\wh v_{2}}{\wh v_1^{\,2}}+\frac{q^2}{1-q^2}}}<1
\ec
\ee
That is to say,
\be
-1+\dis{\frac{1}
{2q\Big(\frac{\wh v_{2}}{\wh v_1^{\,2}}+\frac{q}{1-q^2}\Big)}}<\dis{\frac{\frac{\wh v_{2}}{\wh v_1^{\,2}}+\frac{1}{1-q^2}}{q\frac{\wh v_{2}}{\wh v_1^{\,2}}+\frac{q^2}{1-q^2}}}<1-\dis{\frac{1}
{2q\Big(\frac{\wh v_{2}}{\wh v_1^{\,2}}+\frac{q}{1-q^2}\Big)}}
\ee
or the same
%\be
%-q\Big(\frac{\wh v_{2}}{\wh v_1^{\,2}}+\frac{q}{1-q^2}\Big)+\frac{1}{2}
%<\frac{\wh v_{2}}{\wh v_1^{\,2}}+\frac{1}{1-q^2}<q\Big(\frac{\wh v_{2}}{\wh v_1^{\,2}}+\frac{q}{1-%q^2}\Big)-\frac{1}{2}
%\ee
%This is equivalent to
\be
\frac{\wh v_{2}}{\wh v_1^{\,2}}>\max\left\{
\frac{1/2+1}{q-1},\frac{1/2+\frac{q^2+1}{q^2-1}}{q+1}\right\}=\frac{3}{2}\frac{1}{q-1}
\ee
which is guaranteed by \eqref{eq60} in Corollary\,\ref{cor3.2}.
\ep
We are now ready to complete our proof of the proposition. Indeed,
\be\label{eq85}
\wh\zeta_{X,\F_q;3}^{\,\geq 2}(s)=\frac{1}{2}\wh\zeta_{X,\F_q;3}^{[2]}(s)+ \wh\zeta_{X,\F_q;3}^{[3]}(s)=\wh\zeta_{X,\F_q;3}^{[3]}(s)\left(\frac{\frac{1}{2}\wh\zeta_{X,\F_q;3}^{[3]}(s)}{\wh\zeta_{X,\F_q;3}^{[3]}(s)}+1\right)
\ee
By Lemma\,\ref{lem4.8} just proved, all the zeros of the second factor lie on the line $\Re(s)=\frac{1}{2}-\frac{1}{6}=\frac{1}{3}$. Accordingly, it suffices to show that the zeros of 
$\wh\zeta_{X,\F_q;3}^{[3]}(s)$ cannot be the zeros of $\wh\zeta_{X,\F_q;3}^{[2]}(s)$. 
Recall that $$\ba&  \wh\zeta_{X/\F_q}^{[2]}(s):=\frac{\wh v_{1}^{\,2}}{(1-q^{3s})(1-q^{-3s+3})}\,\wh\zeta_{X,\F_q}(3s-1) \\
& \wh\zeta_{X/\F_q}^{[3]}(s):=
\left( \frac{\wh v_{2}}{(1-q^{-3s+3})}+
\frac{\wh v_{1}^{\,2}}{(1-q^2)(1-q^{-3s+2})}\right)\, \wh\zeta_{X,\F_q}(3s)\\
\ea$$
Obviously, the zeta zeros from the zeta factor $\wh\zeta_{X/\F_q}(3s)$, which are on the line of $\Re(s)=\frac{1}{6}$, cannot be the zeros of $\wh\zeta_{X/\F_q}(3s-1)$, which are on the line $\Re(s)=\frac{1}{2}$, by the Riemann hypothesis for the Artin zeta function $\zeta_{F/\F_q}(s)$. This then leaves  the zeros of the rational function factor in $\wh\zeta_{X,\F_q;3}^{[3]}(s)$, which is clearly not on the line of $\Re(s)=\frac{1}{2}$. Therefore, all zeros of $\frac{1}{2}\wh\zeta_{X,\F_q;3}^{[2]}(s)+ \wh\zeta_{X,\F_q;3}^{[3]}(s)$ are coming from the second factor $\frac{\frac{1}{2}\wh\zeta_{X,\F_q;3}^{[3]}(s)}{\wh\zeta_{X,\F_q;3}^{[3]}(s)}+1$ in \eqref{eq85} and hence lie on the line $\Re(s)=\frac{1}{3}$, as wanted.
\ep

\subsection{Rank three Riemann hypothesis}

Now we are finally ready to complete our proof of Theorem\,\ref{thm4.2} and hence Theorem\,\ref{thm4.3}. We start with the function
$R_{X/\F_q;3}(s):=\left(1+
\frac{
\wh\zeta_{X,\F_q;3}^{\leq 2}(s)}{\wh\zeta_{X,\F_q;3}^{\,\geq 2}(s)}
\right).
$
Recall that
\be
\wh\zeta_{X,\F_q;3}^{\,\geq 2}(s)=\frac{1}{2}
{\wh\zeta_{X,\F_q;3}^{[2]}(s)}
+ \wh\zeta_{X,\F_q;3}^{[3]}(s)=\frac{\prod_{i=1}^{n+2(g-1)}(1-T\ga_i)}{T^{g-1}\prod_{\ell=0}^n(1-q^\ell T)}
\ee
Here, by the Riemann hypothesis for $\wh\zeta_{X,\F_q;3}^{\,\geq 2}(s)$ established in the previous subsection, we have
\be
|\ga_i|=Q^{1/3}\qquad(\forall i=1,\ldots, n+2(g-1))
\ee
Then 
$$\ba
\wh\zeta_{X,\F_q;3}^{\leq 2}(s)=&\wh\zeta_{X,\F_q;3}^{[1]}(s)
+ \frac{1}{2}
{\wh\zeta_{X,\F_q;3}^{[2]}(s)}=\frac{1}{2}
{\wh\zeta_{X,\F_q;3}^{[2]}(1-s)}
+ \wh\zeta_{X,\F_q;3}^{[3]}(1-s)\\
=&\frac{\prod_{i=1}^{n+2(g-1)}(1-\frac{\ga_i}{QT})}{(QT)^{-g+1}\prod_{\ell=0}^n(1-\frac{q^\ell}{QT} )}
%=\frac{\prod_{i=1}^{n+2(g-1)}(1-\frac{\ga_i}{QT})}{(QT)^{-g+1}\prod_{\ell=0}^n(1-\frac{1}{q^{n-\ell}T} )}\\
%=&\frac{\prod_{i=1}^{n+2(g-1)}(1-\frac{\ga_i}{QT})}{(QT)^{-g+1}\prod_{\ell'=n-\ell=n}^0(1-\frac{1}{q^{\ell'}T} )}=\frac{Q^{g-1}q^{1+2+\ldots+n}T^{n+1+g-1}}{(QT)^{n+2(g-1)}}\frac{\prod_{i=1}^{n+2(g-1)}(QT-{\ga_i})}{\prod_{\ell'=n-\ell=n}^0({q^{\ell'}T}-1 )}\\
%=&\frac{Q^{g-1}Q^{(n+1)/2}T^{n+1+g-1}}{(QT)^{n+2(g-1)}}\frac{\prod_{i=1}^{n+2(g-1)}(QT-{\ga_i})}{\prod_{\ell'=n-\ell=n}^0({q^{\ell'}T}-1 )}\\
=\frac{1}{Q^{\frac{n-1}{2}+g-1}T^{g-2}}\frac{\prod_{i=1}^{n+2(g-1)}(QT-{\ga_i})}{\prod_{\ell=0}^n({q^{\ell}T}-1 )}
\ea
$$
Therefore,
$$\ba
\frac{\wh\zeta_{X,\F_q;3}^{\leq[2]}(s)}{\wh\zeta_{X,\F_q;3}^{\,\geq 2}(s)}=&\frac{
{\wh\zeta_{X,\F_q;3}^{[1]}(s)}
+ \frac{1}{2}
{\wh\zeta_{X,\F_q;3}^{[2]}(s)}}
{
\frac{1}{2}
{\wh\zeta_{X,\F_q;3}^{[2]}(s)}
+ \wh\zeta_{X,\F_q;3}^{[3]}(s)}
=\frac{\frac{\prod_{i=1}^{n+2(g-1)}(1-T\ga_i)}{T^{g-1}\prod_{\ell=0}^n(1-q^\ell T)}
}{\frac{1}{Q^{\frac{n-1}{2}+g-1}T^{g-2}}\frac{\prod_{i=1}^{n+2(g-1)}(QT-{\ga_i})}{\prod_{\ell=0}^n({q^{\ell}T}-1 )}}
%=&TQ^{\frac{n-1}{2}+g-1}\frac{{\prod_{i=1}^{n+2(g-1)}(1-T\ga_i)}}
%{{\prod_{i=1}^{n+2(g-1)}(QT-{\ga_i})}}=\frac{T}{\sqrt Q}Q^{\frac{n}{2}+g-1}{\prod_{i=1}^{n+2(g-1)}\frac{1-T\ga_i}{QT-{\ga_i}}}
=\frac{T}{\sqrt Q}{\prod_{i=1}^{n+2(g-1)}\frac{1-T\ga_i}
{\sqrt QT-\frac{\ga_i}{\sqrt Q}}}
\ea$$
We examine the factors $\dis{\frac{1-T\ga_i}
{\sqrt QT-\frac{\ga_i}{\sqrt Q}}}$ under the condition that $|\ga_i|=Q^{1/3}$ for all $1\leq i\leq n+2(g-1)$. Write then $\ga_i=Q^{\frac{1}{3}}e^{i \theta_i}$. Then 
$$\ba
\left|\frac{1-T\ga_i}
{\sqrt QT-\frac{\ga_i}{\sqrt Q}}\right|=\left|\frac{1-TQ^{1/3}e^{i\theta_i}}
{\sqrt QT-\frac{Q^{1/3}e^{-i\theta_i}}{\sqrt Q}}\right|=\left|\frac{1-\sqrt QTQ^{-1/6}e^{i\theta_i}}
{\sqrt QT-{Q^{-1/6}e^{-i\theta_i}}}\right|
\ea$$

Note that
$$\ba
&\left|\frac{1-\sqrt QTQ^{-1/6}e^{i\theta}}
{\sqrt QT-{Q^{-1/6}e^{-i\theta}}}\right|^2
%=&\frac{(1-\sqrt QTQ^{-1/6}e^{i\theta})(1-\sqrt Q\ov TQ^{-1/6}e^{-i\theta})}
%{(\sqrt QT-{Q^{-1/6}e^{-i\theta}})(\sqrt Q\ov T-{Q^{-1/6}e^{i\theta}})}\\
%=&\frac{1-Q^{{1}/{3}}(\ov Te^{-i\theta}+Te^{i\theta})+Q^{2/3}|T|^2}{Q|T|^2-Q^{1/3}(e^{-i\theta} \ov T+e^{i\theta}T)+Q^{-\frac{1}{3}}}=1-\frac{(Q^{{1}/{3}}-1)+(Q-Q^{2/3})|T|^2}{Q|T|^2-Q^{1/3}(e^{-i\theta} \ov T+e^{i\theta}T)+Q^{-\frac{1}{3}}}\\
=1-\frac{(Q^{{1}/{3}}-1)(1+Q^{2/3})|T|^2}{Q|T|^2-Q^{1/3}(e^{-i\theta} \ov T+e^{i\theta}T)+Q^{-\frac{1}{3}}}<1
\ea$$

Obviously, 
\be
\left|\frac{T}{\sqrt Q}\right|\bc<1&\qquad |T|<\sqrt Q\\>1&\qquad |T|>\sqrt Q\ec
\ee
This implies that 
$R_{X/\F_q;3}(s)$ has no zero in the region $\Re(s)<\frac{1}{2}$.  Now note that, by
\S\ref{sec4.1}, we have
\be\label{eq88}
\wh\zeta_{X,\F_q;3}(s)=\wh\zeta_{X/\F_q;3}^{\leq{2}}(s)+\wh\zeta_{X/\F_q;3}^{\geq{2}}(s)=\wh\zeta_{X/\F_q;3}^{\geq2}(s)\left(1+\frac{\wh\zeta_{X/\F_q;3}^{\leq2}(s)}{\wh\zeta_{X/\F_q;3}^{\geq2}(s)}\right)=\wh\zeta_{X/\F_q;3}^{\geq2}(s)\cdot R_{X/\F_q;3}(s)
\ee
Thus, to prove Theorem\,\ref{thm4.2}, what is left to to show that the zeros of $\wh\zeta_{X/\F_q;3}^{\geq2}(s)$ cannot be the zeros of $\wh\zeta_{X,\F_q;3}(s)$. But this is clear, since all zeros of $\wh\zeta_{X/\F_q;3}^{\leq2}(s)$ lie on the line $\Re(s)=1-\frac{1}{3}$ by the functional equation.
\be
\wh\zeta_{X/\F_q;3}^{\leq2}(1-s)=\wh\zeta_{X/\F_q;3}^{\geq2}(s).
\ee
In particular, $\wh\zeta_{X/\F_q;3}^{\leq2}(1-s)$ and $\wh\zeta_{X/\F_q;3}^{\geq2}(s)$ cannot have any common zero. Therefore, from \eqref{eq88},
  the zeros of $\wh\zeta_{X,\F_q;3}(s)$ cannot come from the first  factor $\wh\zeta_{X,\F_q;3}^{\,\geq 2}(s)$, but all come from the second factor $\left(1+\frac{\wh\zeta_{X/\F_q;3}^{\leq2}(s)}{\wh\zeta_{X/\F_q;3}^{\geq2}(s)}\right)$. This proves Theorem\,\ref{thm4.2} and hence also Theorem\,\ref{thm4.3}.

We mention in passing that an analogue of this result holds for the rank $3$ zeta function $\wh\zeta_{\Q;3}(s)$ and hence the $\SL_n$ zeta function $\wh\zeta_{\Q}^{\SL_3}(s)$ the  field $\Q$ of rationals has been proved by Suzuki \cite{Su} (See also \cite{KKS, W2} for general discussions).

\vskip 5.0cm
 Lin WENG
 
 Institute for Fundamental Research
 
 $L$-Academy
 
 \&
  
 Graduate School of Mathematics,
 
 Kyushu University,
 
 Fukuoka 819-0395
 
 Japan
 
 weng@math.kyushu-u.ac.jp

\end{document}